\documentclass[12pt,letterpaper]{amsart}
\usepackage{fontenc,amsmath,mathrsfs, amsthm, verbatim,amssymb, MnSymbol, enumerate, multicol, hyperref, geroENG, color}
\usepackage[english]{babel}
\allowdisplaybreaks
\usepackage{times}
\allowdisplaybreaks
\usepackage[top=2.5cm, bottom=2.5cm, left=2.25cm, right=2.25cm]{geometry}

\begin{document}
\title[Measure theoretic properties of large products]{Measure theoretic properties of large products of consecutive partial quotients}
\author[Brown-Sarre]{Adam Brown-Sarre}
\author[Gonz\'alez Robert]{Gerardo Gonz\'alez Robert$^\ast$}
\author[Hussain]{Mumtaz Hussain}
\address[ABS, GGR, MH]{Department of Mathematical and Physical Sciences,  La Trobe University, Bendigo 3552, Australia. }
\address[GGR]{Facultad de Ciencias, Universidad Nacional Autónoma de México, Investigación Científica, C.U., Coyoacán, 04510 Ciudad de México, Mexico}
\email[ABS]{20356213@students.latrobe.edu.au}
\email[GGR]{g.robert@latrobe.edu.au, gero@ciencias.unam.mx}
\email[MH]{m.hussain@latrobe.edu.au}
\thanks{$^\ast$Corresponding author.}


\begin{abstract}
The theory of uniform approximation of real numbers motivates the study of products of consecutive partial quotients in regular continued fractions. For any non-decreasing positive function $\varphi:\Na\to\RE_{\geq 2}$, we determine the Lebesgue measure of the set $\mathcal{F}_{\ell}(\vphi)$ of irrational numbers $x$ whose regular continued fraction $x~=~[a_1(x),a_2(x),\ldots]$ is such that, for infinitely many $n\in\Na$, there are two numbers $1\leq j<k \leq n$ satisfying
\[
a_{k}(x)\cdots a_{k+\ell-1}(x) \geq \varphi(n), \;
a_{j}(x)\cdots a_{j+\ell-1}(x) \geq \varphi(n). 
\]
This result generalizes previous work by Tan and Zhou (Nonlinearity, 2024).
A consequence of our result is that the strong law of large numbers for products of $\ell$ consecutive partial quotients is impossible even if the block with the largest product is removed.
We also compute the Hausdorff dimension of $\clF_3(\vphi)$.

\end{abstract}
\maketitle
\section{Introduction}

Regular continued fractions are a representation system of real numbers with outstanding approximation properties. 
It is well-known that, for each $x\in \RE\setminus\QU$ there is a single sequence of integers $(a_j(x))_{j\geq 0}$ with $a_j(x)\geq 1$ for $j\geq 1$ such that
\[
x
=
a_0(x) + \cfrac{1}{a_1(x) + \cfrac{1}{a_2(x) + \cfrac{1}{\ddots}}}
=
[a_0(x); a_1(x),a_2(x),\ldots].
\]
The numbers $a_n(x)$ are the \textbf{partial quotients} of $x$. 
A necessary and sufficient condition for $x\in [0,1)$ is $a_0(x)=0$; in this case, we denote the continued fraction of $x$ by $[a_1(x),a_2(x),\ldots]$. 
A classical result (see, for example, \cite[Equation (34)]{KhinchinCFBook}) says that the $n$th \textbf{convergent} of $x$, given by
\[
\frac{p_n(x)}{q_n(x)}
=
[a_1(x),\ldots, a_n(x)]
=
\cfrac{1}{a_1(x)+\cfrac{1}{\ddots + \cfrac{1}{a_n(x)}}},
\]
where $q_n(x)\geq 1$ and $p_n(x)$ are coprime integers, satisfies 
\[
\frac{1}{(a_{n+1}(x)+2)q_n^2(x)}
<
\left| x - \frac{p_n(x)}{q_n(x)}\right|
<
\frac{1}{a_{n+1}(x)q_n^2(x)}.
\]
Then, since the convergents of $x$ are its best rational approximations, the rate at which rationals approximate $x$ is closely related to the size of its partial quotients. 
This leads us to consider, for any function $\vphi:\Na\to\RE_{>0}$, the set
\[
\clE_1(\varphi)
:= \{x \in [0,1):a_n(x) \geq \varphi(n) \text{ for i.m. } n\in \mathbb{N}\}.
\]
Here and in what follows, ``i.m.'' stands for ``infinitely many,'' also, $\leb$ stands for the Lebesgue measure on $[0,1]$.  
The Borel-Bernstein Theorem establishes that
\[
\leb\left( \clE_1(\vphi) \right)
=
\begin{cases}
    0, \text{ if }  \displaystyle\sum_{n\geq 1} \frac{1}{\vphi(n)}<\infty,\\
    1, \text{ if }  \displaystyle\sum_{n\geq 1} \frac{1}{\vphi(n)}=\infty,\\
\end{cases}
\]
(see, for example, \cite[Theorem 1.11]{BugeaudBookCUP160}). 
{\L}uczak \cite{Luc1997} and Feng et al. \cite{FengWuLiang1997} calculated the Hausdorff dimension, denoted as $\hdim$, of $\clE_1(\vphi)$ when $\vphi(n)=c^{b^n}$ for some fixed $c,b>1$ (see Lemma \ref{Lem:Luczak}). 
In \cite{WangWu2008}, Wang and Wu computed $\hdim \clE_1(\vphi)$ for a general $\vphi$ (Theorem \ref{TEO:HWX:HD} below with $\ell=1$).
While studying improvements to the classical Dirichlet theorem, Kleinbock and Wadleigh \cite{KleinbockWadleigh2019} showed that the set of Dirichlet non-improvable numbers is related to the set
\[
\clE_{2}(\varphi)
:= \{x \in [0,1):a_n(x)a_{n+1}(x) \geq \varphi(n) \text{ for i.m.   } n\in \mathbb{N}\}
\]
and computed $\leb(\clE_{2}(\varphi))$. 
In \cite{HuangWuXu2020}, Huang, Wu, and Xu extended this result to blocks of arbitrary length; namely, for $\vphi:\Na\to\RE_{>0}$ and $\ell\in \Na$, they obtained the Lebesgue measure and Hausdorff dimension of
\[
\clE_{\ell}(\varphi)
:= \{x \in [0,1):a_n(x)\cdots a_{n+\ell-1}(x) \geq \varphi(n) \text{ for i.m.   } n\in \mathbb{N}\}.
\]
\begin{teo01}[Huang, Wu, Xu \cite{HuangWuXu2020}]\label{TEO:HWX:MEAS}
If $\vphi:\Na\to \RE_{\geq 2}$ and $\ell\in\Na$, then
\[
\leb\left( \clE_{\ell}(\vphi)\right)
=
\begin{cases}
0, &\text{\rm if } \displaystyle\sum_{n=1}^{\infty} \frac{\log^{\ell-1}\vphi(n)}{\vphi(n)}<\infty, \\    
1, &\text{\rm if } \displaystyle\sum_{n=1}^{\infty} \frac{\log^{\ell-1}\vphi(n)}{\vphi(n)}=\infty.
\end{cases}
\]
\end{teo01}
To state the Hausdorff dimension result, for any $\vphi:\Na\to \RE_{\geq 2}$, let $B_{\vphi}, b_{\vphi}\in [1, \infty]$ be determined by
\begin{equation}\label{Bb}
\log B_{\vphi} := \liminf_{n \xrightarrow{} \infty} \frac{\log \varphi(n)}{n} 
\;\text{ and }\;
\log b_{\vphi} = \liminf_{n \xrightarrow{} \infty} \frac{\log \log \varphi(n)}{n}.
\end{equation}
For the rest of the paper, $T$ stands for the Gauss map (see Subsection \ref{Subsec:ContFrac}).
The pressure functions $P(T,\cdot)$ appearing in Theorem \ref{TEO:HWX:HD} and Theorem \ref{TEO:MAIN:HD} are defined in Subsection \ref{Subsec:PressureFunctions}.
\begin{teo01}[Huang, Wu, Xu \cite{HuangWuXu2020}]\label{TEO:HWX:HD}
Let $(f_n)_{n\geq 1}$ be the sequence of real functions on $[0,1]$ given by
\[
f_1(s)=s, 
\;\text{ and }\;
f_{n+1}=\frac{sf_n(s)}{1-s+f_n(s)} \;\text{ for all } n\in\Na.
\]
For any $\vphi:\Na\to\RE_{\geq 2}$ and $\ell\in\Na$, then the Hausdorff dimension of $\clE_{\ell}(\vphi)$ is 
\[
\hdim \clE_{\ell}(\vphi)
=
\begin{cases}
1, & \text{\rm if } B_{\vphi}=1, \\
\inf\{s \geq 0: P(T,-f_{\ell}(s)\log B_{\vphi} - s \log |T'|) \leq 0\},  & \text{\rm if }  1<B_{\vphi}<\infty,\\
\frac{1}{1+b_{\vphi}},  & \text{\rm if }  B_{\vphi}=\infty.
\end{cases}
\]
\end{teo01}

Let us adopt a probabilistic perspective and interpret partial quotients as random variables on $[0,1)$; for our purposes, their definition on rational numbers is irrelevant. 
These random variables are not independent; hence, the Borel-Bernstein Theorem is not a trivial consequence of the Borel-Cantelli Lemma (Lemma \ref{Le:BorelCantelli}). 
Moreover, the Borel-Bernstein theorem tells us that the inequality $a_n(x)\geq n\log n$ holds infinitely often for almost every number $x$ with respect to $\leb$. This implies that, if $S_n(x):=a_1(x) + \cdots + a_n(x)$ for $n\in\Na$, we have
\[
\limsup_{n\to\infty} \frac{1}{n} S_n(x) = \infty
\;\text{ for almost all } x\in [0,1).
\]
Despite not being independent and having infinite expected value, i.e., $\int a_n\md\leb =\infty$ for all $n$, the partial quotients satisfy a weak law of large numbers: Khinchin \cite[\S. 4]{Khintchine1935} proved that $S_n/(n\log n)$ converges in Lebesgue measure to $\frac{1}{\log 2}$ as $n\to\infty$. Furthermore, Philipp \cite[Theorem 1]{Philipp1988} proved that there is no reasonable function $\sigma$ for which $S_n/\sigma(n)$ converges almost everywhere to a finite positive constant. Large partial quotients are to blame for this phenomenon; in fact, Diamond and Vaaler \cite[Corollary 1]{DiamondVaaler1986} showed that convergence almost everywhere does hold for Khinchin's averaging function when the largest partial quotient is omitted; that is, 
\begin{equation*}
\lim_{n\to\infty} 
\frac{1}{n\log n}\left( S_n(x) - \max_{1\leq j \leq n} a_j(x)\right)
=
\frac{1}{\log 2}
\;\text{ for $\leb$-almost all } x\in [0,1).
\end{equation*}
This result was extended by Hu, Hussain, and Yu for sums of products of two partial quotients \cite[Theorem 1.5]{HuHussainYu2021}. They showed that almost every $x\in [0,1)$ satisfies
\begin{equation*}
\lim_{n\to\infty} 
\frac{1}{n\log^2 n}\left( \sum_{j=1}^n a_j(x)a_{j+1}(x) - \max_{1\leq j \leq n} a_j(x)a_{j+1}(x)\right)
=
\frac{1}{2\log 2}.
\end{equation*}

A key aspect in the proofs of \eqref{Eq:Thm:DV} and \eqref{Eq:Thm:HHY}, respectively, is showing that, for $c>\frac{1}{2}$, we have
\begin{equation}\label{Eq:Thm:DV}
\leb\left(
\left\{
x\in [0,1):
\begin{matrix}
a_j(x) \geq  n (\log n)^{c}, \\
a_k(x) \geq  n (\log n)^{c},
\end{matrix}\;
\text{ for some } 1\leq j < k \leq n \text{ for i.m.   } n\in\Na
\right\}\right)=0
\end{equation}
and, for $d>\frac{3}{2}$,
\begin{equation}\label{Eq:Thm:HHY}
\leb\left(
\left\{x\in [0,1): 
\begin{matrix}
a_j(x)a_{j+1}(x) \geq n (\log n)^{d} , \\
a_k(x)a_{k+1}(x) \geq n (\log n)^{d},  \\
\end{matrix}\;
\text{ for some } 1\leq j < k \leq n \text{ for i.m.   } n\in\Na\right\}\right)=0.
\end{equation}
Recently, Tan, Tian, and Wang \cite{TanTianWang2023} generalized \eqref{Eq:Thm:DV} to a wider family of functions (Theorem \ref{TEO:TTW:MEAS}), and Tan and Zhou generalized \eqref{Eq:Thm:HHY} in \cite{TanZhou2024} (Theorem \ref{TEO:TZ:MEAS}). 
To be precise, given a non-decreasing function $\vphi:\Na\to\RE_{\geq 2}$ and $\ell\in\Na$, define 
\[
\clF_{\ell}(\vphi)
:=
\left\{ x\in [0,1):
\begin{matrix}
a_{j}(x)\cdots a_{j+\ell-1}(x) \geq \vphi(n), \\
a_{k}(x)\cdots a_{k+\ell-1}(x) \geq \vphi(n),
\end{matrix}\;
\text{ for some } 1\leq j < k \leq n \text{ for i.m.   } n\in\Na\right\}.
\]
The assumption $\vphi(n)\geq 2$ for $n\in\Na$ is not restrictive at all, since  $\clF_{\ell}(\vphi)=[0,1)\setminus \QU$ whenever $\vphi(n)\leq 1$ for infinitely many $n\in \Na$. 

\begin{teo01}[Tan, Tian, Wang \cite{TanTianWang2023}]\label{TEO:TTW:MEAS} 
If  $\vphi:\Na\to[2,\infty)$ is non-decreasing, then 
    \[
    \leb(\clF_1(\vphi))
    =
    \begin{cases}
        0, &\text{\rm if }  \displaystyle\sum_{n\geq 1} \frac{n}{\varphi(n)^2}<\infty,\\
        1, &\text{\rm if }  \displaystyle\sum_{n\geq 1} \frac{n}{\varphi(n)^2}=\infty.
    \end{cases}
    \]
\end{teo01}

\begin{teo01}[Tan, Zhou, \cite{TanZhou2024}]\label{TEO:TZ:MEAS}
If  $\vphi:\Na\to[2,\infty)$ is non-decreasing, then 
\[
\leb(\clF_2(\vphi))
=
\begin{cases}
0, &\text{\rm if } \displaystyle\sum_{n\geq 1} 
\left(\frac{n\log^2  \vphi(n)}{\vphi^2(n)} + \frac{1}{\vphi(n)}\right) <\infty, \\   
1, &\text{\rm if } \displaystyle\sum_{n\geq 1} 
\left(\frac{n\log^2  \vphi(n)}{\vphi^2(n)} + \frac{1}{\vphi(n)}\right) =\infty.
\end{cases}
\]
\end{teo01}
Our first result gives the Lebesgue measure of $\clF_{\ell}(\vphi)$ for any $\ell\in\Na_{\geq 2}$.
\begin{teo01}\label{TEO:MAIN:MEASURE}
If  $\vphi:\Na\to[2,\infty)$ is non-decreasing and $\ell\in\Na_{\geq 2}$, then
\[
\leb(\clF_{\ell}(\vphi))
=
\begin{cases}
0 &\text{\rm if } 
\displaystyle\sum_{n\geq 1} 
\left(\frac{n\log^{2(\ell-1)}\vphi(n)}{\vphi^2(n)} + \frac{\log^{\ell-2} \vphi(n)}{\vphi(n)}\right)<\infty, \\
1 &\text{\rm if } \displaystyle\sum_{n\geq 1} 
\left(\frac{n\log^{2(\ell-1)}\vphi(n)}{\vphi^2(n)} + \frac{\log^{\ell-2} \vphi(n)}{\vphi(n)}\right)=\infty.
\end{cases}
\]
\end{teo01}
\begin{rema01}
To prove Theorem \ref{TEO:MAIN:MEASURE}, we need new estimates corresponding to overlapping blocks with large products. 
The lengths of such overlapping blocks range from $0$ up to $\ell-1$, whereas in Theorem \ref{TEO:TZ:MEAS} they are only $0$ or $1$, and they do not exist in Theorem \ref{TEO:TTW:MEAS}. 
The biggest novelty in our proof is Theorem \ref{Thm:OverlapEstimatesSeries}.
\end{rema01}

\begin{rema01}
It is well-known that the measure $\mu$ given by $\md \mu = \frac{1}{\log 2} \frac{1}{1+x} \md \leb$ is $T$-ergodic (see, for example, \cite[Theorem 3.7]{EinsiedlerWardET}).
Then, Birkhoff's Ergodic Theorem yields
\[
\lim_{n\to\infty}
\frac{1}{n}
\sum_{j=1}^n a_j(x)\cdots a_{j+\ell-1}(x)
=
\infty
\;
\text{ for almost every }x\in [0,1)
\]
(see, for example, \cite[Corollary 3.8]{EinsiedlerWardET}).
Moreover, by Theorem \ref{TEO:MAIN:MEASURE}, for almost every $x\in [0,1)$ there are infinitely many numbers $n\in\Na$ for which there exist $j,k\in\Na$, $1\leq j<k\leq n$ satisfying 
\[
a_{j}(x)\cdots a_{j+\ell-1}(x)  \geq n\log^{2}n
\;\text{ and }\;
a_{k}(x)\cdots a_{k+\ell-1}(x)  \geq n\log^{2}n;
\]
thus, for such $x$,
\[
\limsup_{n\to\infty}  
\frac{1}{n} \left(\sum_{j=1}^n a_{j}(x)\cdots a_{j+\ell-1}(x)  
- 
\max_{1\leq k\leq n} a_{k}(x)\cdots  a_{k+\ell-1}(x) \right)
=
\infty.
\]
We suspect that Theorem \ref{TEO:MAIN:MEASURE} is a crucial step toward the extension of \eqref{Eq:Thm:DV} and \eqref{Eq:Thm:HHY}; in fact, we believe that it will allow us to show that 
\[
\lim_{n\to\infty} 
\frac{1}{n\log^{\ell} n} 
\left( \sum_{j=1}^na_{j}(x)\cdots a_{j+\ell-1 }(x) 
- 
\max_{1\leq k\leq n} a_{k}(x)\cdots  a_{k+\ell-1}(x)\right)
=
\frac{1}{\ell \log 2}
\;\text{ for almost all } x\in [0,1).
\]
We state other anticipated results in Section \ref{Section:FinalRemarks}. 

\end{rema01}

In Diophantine approximation, the Hausdorff dimension is regularly used to distinguish among null sets.
In \cite{TanTianWang2023}, the authors calculated $\hdim \clF_{1}(\vphi)$ (Theorem \ref{TEO:TTW:HD}) and $\hdim \clF_{2}(\vphi)$ (Theorem \ref{TEO:TZ:HD}) was calculated in \cite{TanZhou2024}. 

\begin{teo01}[Tan, Tian, Wang \cite{TanTianWang2023}]\label{TEO:TTW:HD}
If $\vphi:\Na\to [2,\infty)$ is non-decreasing, then
\[
\hdim \clF_1(\vphi)  
=
\begin{cases}
1, &\text{\rm if } B_{\vphi}=1, \\
\inf\{s \geq 0: P(T,-(3s-1)\log B_{\vphi} - s \log |T'|) \leq 0\}, &\text{\rm if }  1<B_{\vphi}<\infty,\\
\frac{1}{1+b_{\vphi}}, &\text{\rm if }  B_{\vphi}=\infty.
\end{cases}  
\]
\end{teo01}

\begin{teo01}[Tan, Zhou \cite{TanZhou2024}]\label{TEO:TZ:HD}
If $\vphi:\Na\to [2,\infty)$ is non-decreasing, then
\[
\hdim \clF_2(\vphi)  
=
\begin{cases}
1, &\text{\rm if }  B_{\vphi}=1, \\
\inf\{s \geq 0: P(T,-(3s-1-s^2)\log B_{\vphi} - s \log |T'|) \leq 0\}, &\text{\rm if }  1<B_{\vphi}<\infty,\\
\frac{1}{1+b_{\vphi}}, &\text{\rm if } B_{\vphi}=\infty.
\end{cases}  
\]
\end{teo01}
Our second result, Theorem \ref{TEO:MAIN:HD}, is the computation of $\hdim \clF_3(\vphi)$. 
Define $g:[0, 1]\to\RE$ by
\[
g(s)
:=
\frac{3s^3-5s^2+4s-1}{s^2-s+1}.
\]

\begin{teo01}\label{TEO:MAIN:HD}
If $\vphi:\Na\to [2,\infty)$ is non-decreasing, then
\[
\hdim \clF_3(\vphi)  
=
\begin{cases}
1, &\text{\rm if }  B_{\vphi}=1, \\
\inf\{s \geq 0: P(T,-g(s)\log B_{\vphi} - s \log |T'|) \leq 0\}, &\text{\rm if }  1<B_{\vphi}<\infty,\\
\frac{1}{1+b_{\vphi}}, &\text{\rm if }  B_{\vphi}=\infty.
\end{cases}  
\]
\end{teo01}

The paper is organized as follows. Section \ref{Section:Preleminaries} contains some preliminaries and auxiliary lemmas from measure theory and continued fractions, in particular on the Hausdorff dimension and pressure functions. In Section \ref{Section:Series}, we estimate certain series that are fundamental for our proofs. The proof of Theorem \ref{TEO:MAIN:MEASURE} is in Section \ref{Section:ProofMeasure} and the proof of Theorem \ref{TEO:MAIN:HD} is in Section \ref{Section:ProofHD}. Lastly, Section \ref{Section:FinalRemarks} contains some final remarks and directions to further this study.  

\paragraph{Notation.} 
Whenever $(x_n)_{n\geq 1}$ and $(y_n)_{n\geq 1}$ are two sequences of positive real numbers, we write $x_n\ll y_n$ to mean that for some constant $c>0$ and every $n\in\Na$ we have $x_n\leq cy_n$. 
When $c$ depends on some parameter $\lambda$, we write $x_n\ll_{\lambda} y_n$. By $x_n\asymp y_n$ we mean $x_n\ll y_n$ and $y_n\ll x_n$. 
We write $x_n\asymp y_n$ when $x_n\ll y_n$ and $y_n\ll x_n$ and we write $\asymp_{\lambda}$ when the constants depend on some parameter $\lambda$.
For the rest of the paper, $\ell$ will denote an arbitrary natural number.

\medskip

\noindent{\bf Acknowledgements.} 
GGR and MH were funded by the Australian Research Council Discovery Project Grant Number 200100994. 
GGR was funded by the Secretaría de Ciencia, Humanidades, Tecnología e Innovación (Mexico) through the Sistema Nacional de Investigadores e Investigadoras. 
We thank Rafael Alcaraz Barrera, Nikita Shulga, and Ben Ward, and especially Valeriia Starichkova for their important observations.
\section{Preliminaries}\label{Section:Preleminaries}
In this section, we recall some definitions and elementary results on measure theory, Hausdorff dimension, continued fractions, and pressure functions. 

\subsection{Full and null sets}
Recall that the limsup set of a sequence $(E_n)_{n\geq 1}$ of subsets of a given set $X$ is
\[
\limsup_{n\to\infty} E_n = \bigcap_{n\geq 1} \bigcup_{m\geq n} E_m.
\]

\begin{teo01}[Borel-Cantelli Lemma]\label{Le:BorelCantelli}
Let $(X, \scA, \mu)$ be a probability space and $(E_k)_{k\geq 1}$ a sequence of measurable sets.
\begin{enumerate}[\rm i.]
\item If $\displaystyle\sum_{n\geq 1}\mu(E_n)<\infty$, then $\mu(\displaystyle\limsup_{n\to\infty} E_n)=0$.
\item If the sets $E_n$, $n\in\Na$, are independent and $\displaystyle\sum_{n\geq 1}\mu(E_n)=\infty$, then $\mu(\displaystyle\limsup_{n\to\infty} E_n)=1$.
\end{enumerate}
\end{teo01}
\begin{proof}
  See, for example, \cite[Lemma 4.18]{KallenbergBook2021}.
\end{proof}

The independence assumption is rather restrictive; it will not hold in our context. 
A way to circumvent this problem is to show the positive measure of the limsup set and that the measure is either full or null by some other means. 
We obtain the full measure statement in Theorem \ref{TEO:MAIN:MEASURE} by showing that our limsup set intersects a large portion of every interval and applying Knopp's Lemma (Lemma \ref{Le:Knopp}). 

\begin{lem01}[Knopp's Lemma]\label{Le:Knopp}
Let $A \subseteq [0,1)$ be a Borel set and $\mathcal{C}$ a class of subintervals of $[0, 1)$ satisfying the following conditions:
    \begin{enumerate}
        \item Every open subinterval of $[0,1)$ is an at most countable union of disjoint elements from $\mathcal{C}$,
        \item There is a constant $c>0$ such that for any $B \in \mathcal{C}$ we have $\leb (A\cap B) \geq c \leb(B)$.
    \end{enumerate}
    Then, $\leb(A)=1$.
\end{lem01}
\begin{proof}
    See, for example, \cite[Lemma 2.3.1]{MR4701087}
\end{proof}

To apply Lemma \ref{Le:Knopp}, we need the following result due to Chung and Erd{\H o}s \cite{ChungErdos1952} (see also \cite[Lemma 6]{Sprindzuk1979} for the formulation below).

\begin{lem01}[Chung, Erd{\H o}s, \cite{ChungErdos1952}]\label{Le:ErdosChung}
Let $(X, \scA, \mu)$ be a probability space and $(E_k)_{k\geq 1}$ a sequence of measurable sets. 
If $\displaystyle\sum_{n\geq 1}\mu(E_n)=\infty$, then
\[
\mu\left(\limsup_{n\to\infty} E_n\right)
\geq 
\limsup_{N\to\infty} \frac{\left( \displaystyle\sum_{1\leq n\leq N} \mu(E_n)\right)^2}{\displaystyle\sum_{1\leq i, j\leq N} \mu(E_i\cap E_j)}.
\]
\end{lem01}

\subsection{Hausdorff dimension}
We recall the definition of the Hausdorff dimension and establish some notation.
For the general theory of the subject, see \cite{FalconerBookFractal}.
Given a bounded subset $A$ of $\RE$, let $|A|$ be its diameter, that is, 
\[
|A|
=
\begin{cases}
\sup\{|a-b|:a,b\in A\}, &\text{ if } A\neq \vac, \\
0, &\text{ if } A= \vac.
\end{cases}
\]
For any $\delta > 0$, a \textbf{$\delta$-cover} of $A$ is a countable collection $\mathfrak{B}$ of bounded subsets of real numbers such that 
\[
\sup_{B\in \mathfrak{B}} |B|\leq \delta
\quad\text{ and }\quad
A\subseteq \bigcup_{ B\in \mathfrak{B} } B.
\]
In this case, for $s>0$, we define
\[
\Lambda_s(\mathfrak{B}) = \sum_{B\in \mathfrak{B}} |B|^s
\]
and 
\[
\clH^s_{\delta}(A)
=
\inf\left\{\Lambda_{s}(\mathfrak{B}) : \mathfrak{B} \text{ is a }\delta-\text{cover of } A\right\}.
\]
The \textbf{$s$-Hausdorff measure} $\clH^s$ of $A$, denoted by $\clH^s(A)$, is 
\[
\clH^s(A)
:=
\lim_{\delta\to 0} \clH^s_{\delta}(A).
\]
It is well-known that there is a number, called the \textbf{Hausdorff dimension of $A$} and denoted by $\hdim A$, such that 
\[
\clH^s(A)
:=
\begin{cases}
\infty, & \text{if} \ \ s< \hdim A,\\
0, & \text{if} \ \  s> \hdim A.
\end{cases}
\]

\subsection{Continued fractions}\label{Subsec:ContFrac}
For each real number $t$, let $[t]$ be its integer part. 
Regular continued fractions can be generated with the Gauss map $T:[0,1)\to [0,1)$ given by 
\[
T(x) = 
\begin{cases}
x^{-1} - [x^{-1}], & x\neq 0, \\
0, & x=  0.
\end{cases}
\]
For any $x\in [0,1)\setminus \QU$, let $a_1:[0,1)\setminus \QU\to\Na$ be given by $a_1(x)= [x^{-1}]$ and define $a_n(x)=a_1(T^{n-1}(x))$ for each $n\in\Na$. 
As in the introduction, the natural numbers $(a_n(x))_{n\geq 1}$ are the \textbf{partial quotients} of $x$. 
This is the only sequence of natural numbers satisfying
\[
x=
[a_1(x),a_2(x),\ldots]
=
\lim_{n\to\infty} [a_1(x),a_2(x),\ldots, a_n].
\]
For any non-negative integer $n$, the $n$th \textbf{convergent} of $x$ is 
\[
\frac{p_n(x)}{q_n(x)}
:=
[a_1(x), \ldots, a_n(x)].
\]
An elementary result (see, for example, \cite[Proposition 1.1]{HensleyBook}) states that the sequences $(p_n(x))_{n\geq 0}$, $(q_n(x))_{n\geq 1}$ are given by
\[
q_0(x)=1, \quad
q_{1}(x)=a_1(x), \quad
p_0(x)=0, \quad
p_1(x)=1,
\]
and 
\[
p_{n+1}(x)=a_{n+1}(x)p_{n}(x) + p_{n-1}(x), \quad
q_{n+1}(x)=a_{n+1}(x)q_{n}(x)+q_{n-1}(x)
\;\text{ for all } n\in\Na.
\]
The terms of the sequence $(q_n(x))_{n\geq 1}$ are called \textbf{continuants}. 
Clearly, the $n$th continuant depends only on the first $n$ partial quotients of $x$. Therefore, any other number $y$ whose first $n$ partial quotients are exactly those of $x$ has the same first $n$ continuants as $x$. 
For each $n\in\Na$ and $\bfa=(a_1,\ldots, a_n)\in\Na^n$, the \textbf{fundamental interval} $I_n(\bfa)$ is
\[
I_n(\bfa)
:=
\left\{ x\in [0,1): a_1(x)=a_1,\ldots, a_n(x)=a_n\right\}.
\]
In the following, we denote by $q_n(\bfa)$ the $n$-th continuant of any number in $I_n(\bfa)$. We omit the argument of $q_n$ if there is no risk of ambiguity.
\begin{lem01}[{\cite[Section 12]{KhinchinCFBook}}]\label{Lem:CF01}
There are absolute constants such that for all $n\in\Na$ and all $\bfa\in\Na^n$ we have
\[
|I_n(\bfa)|
\asymp
\frac{1}{q_n^2(\bfa)}.
\]
\end{lem01}
For any $n,m\in\Na$, $\bfa=(a_1,\ldots, a_n)\in\Na^n$ and $\bfb=(b_1,\ldots, b_m)\in\Na^m$, we write
\[
(\bfa,\bfb)
:=
(a_1,\ldots, a_n,b_1,\ldots, b_m).
\]
\begin{lem01}[{\cite[Proposition 1.1]{HensleyBook}}]\label{Lem:CF02}
For any $m,n\in\Na$,  $\bfa\in\Na^n$, and $\bfb\in\Na^m$, 
\[
q_{n}(\bfa)q_{m}(\bfb)
\leq
q_{n+m}(\bfa,\bfb)
\leq 
2q_{n}(\bfa)q_{m}(\bfb).
\]
\end{lem01}
\begin{lem01}[{\cite[Lemma 2.1]{Wu2006}}]\label{Lem:CF03}
For any $n\in\Na$,  $\bfa=(a_1,\ldots, a_n)\in\Na^n$, and $k\in\{1,\ldots,n\}$, we have
\[
\frac{a_{k} + 1}{2}
\leq 
\frac{q_n(a_1,\ldots, a_n)}{q_{n-1}(a_1,\ldots,a_{k-1},a_{k+1},\ldots, a_n)}
\leq 
a_{k}+1.
\]
\end{lem01}

\begin{lem01}[\cite{FengWuLiang1997, Luc1997}]\label{Lem:Luczak}
For any constants $b,c>1$, we have 
\begin{align*}
\frac{1}{b+1}
&=
\hdim
\left\{ x\in [0,1): a_n(x)\geq b^{c^n} \text{ for all large } n\in\Na \right\}, \\
&=
\hdim
\left\{ x\in [0,1): a_n(x)\geq b^{c^n} \text{ for i.m.   } n\in\Na\right\}.
\end{align*}
\end{lem01}
Although not explicitly stated, the next classical result can be retrieved from the discussion in \cite[Chapter 2]{KhinchinCFBook}.
\begin{lem01}\label{PR:UnionFundamentalIntervals}
For $n\in\mathbb{N}$, $\bfa \in \mathbb{N}^{n}$, and $a,b\in\mathbb{N}$, we have
\[
\leb\left(
\bigcup_{a\leq j \leq b} I_{n+1}(\bfa,j)
\right)
=
\frac{(b+1)-a}{(aq_n +q_{n-1})((b+1)q_n +q_{n-1})}.
\]
\end{lem01}

\subsection{Pressure functions}\label{Subsec:PressureFunctions}
Pressure functions have become popular in Diophantine approximation because their zeros are the Hausdorff dimension of certain sets determined by dynamical means. In this paper, we restrict ourselves to the dynamical system $([0,1)\setminus \QU,T)$. 
For the general theory of pressure functions, see the monograph by Mauldin and Urba{\'n}ski \cite{MauldinUrbanskiBook}. 

Given a function $\psi:[0,1] \to\RE$, for any $x\in [0,1]$ and $n\in\Na$, we let $S_n(\psi;x)$ be the ergodic sum
\[
S_n(\psi;x)
:=
\psi(x) + \psi(T(x)) + \ldots + \psi\left(T^{n-1}(x)\right).
\]
The $n$th variation of $\psi$ is given by
\[
\Var_n(\psi)
:=
\sup\left\{ |\psi(x)-\psi(y)|:I_n(x)=I_n(y)\right\}.
\]
For any non-empty subset $A\subseteq \Na$, we consider
\[
X_A
:=
\left\{ x\in [0,1): a_n(x)\in A \text{ for all } n\in\Na\right\}.
\]
The pressure function $P_{A}(T,\cdot)$ on $\psi$ is given by
\[
P_{A}(T,\psi)
:=
\lim_{n\to\infty} 
\frac{1}{n} 
\log\left( 
\sum_{\bfa\in A^n} \sup_{x\in I_n(\bfa)\cap X_{ A}} \exp\left( S_n(\psi,x)\right)
\right).
\]
In this context, we call $\psi$ a \textbf{potential}. 
When $A=\Na$, we write $P(T,\psi)$ instead of $P_{\Na}(T,\psi)$. The existence of $P_{A}(T,\psi)$ is not immediate, but it does exist under some regularity conditions. Furthermore, we can approximate $P(T,\psi)$ by $P_{A}(T,\psi)$ with $A\subseteq \Na$ finite.

\begin{lem01}[{\cite[Proposition 2.4]{LiWangWuXu2014}}]\label{Lem:PF:01}
Let $\psi:[0,1]\to\RE$ be such that $\Var_1(\psi)<\infty$ and $\Var_n(\psi)\to 0$ as $n\to\infty$. Then, $P_{A}(T,\psi)$ is well defined.
Moreover, the supremum in the definition of $P_{A}(T,\psi)$ may be replaced with $\exp(S_n(\psi;x))$ for any $x\in I_n(\bfa)\cap X_{A}$.
\end{lem01}
\begin{lem01}[{\cite[Proposition 2]{HanusMauldinUrbanski2002}}]\label{Lem:PF:02}
If $\psi:[0,1]\to\RE$ is such that $\Var_1(\psi)<\infty$ and $\Var_n(\psi)\to 0$ as $n\to\infty$, then
\[
P(T,\psi)
=
\sup \left\{ P_{A}(T,\psi): A\subseteq \Na\;\text{\rm is finite}\right\}.
\]
\end{lem01}
For any continuous function $f:[0,1]\to\RE$ and $0<s<1$, consider the potential
\[
\psi_s(x)
:= 
-f(s)\log B - s\log|T'(x)|.
\]
Then, $\psi_s$ satisfies the hypotheses of lemmas \ref{Lem:PF:01} and \ref{Lem:PF:02}. Since for every $x\in [0,1)$ with $n$ partial quotients we have
\[
(T^n)'(x)
=
\frac{(-1)^n}{(xq_{n-1}-p_{n-1})^2},
\]
the pressure function $P(T,\cdot)$ evaluated on $\psi_s$ becomes
\[
P(T,\psi_s)
=
\lim_{n\to\infty} 
\frac{1}{n} 
\log\left( 
\sum_{\bfa\in A^n} \frac{1}{B^{nf(s)} q_{n}^{2s}(\bfa)} 
\right).
\]
In this paper, we work with potentials $\psi_s$ given by 
\[
\psi_s(x)
=
-g(s)\log B -s\log |T'(x)|,
\]
where $B>1$ is a given real number, $0<s<1$, and $g$ is as in Theorem \ref{TEO:MAIN:HD}.
Define
\[
s_B
:=
\inf\{s\geq 0: P(T,-g(s)\log B -s\log |T'(x)|)\leq 0\}
\]
and, for $m\in\Na$,
\begin{equation}\label{Eq:DefSmB}
s_{m,B}
:=
\inf
\left\{ s\geq 0:
\sum_{\bfa\in\Na^m} 
\left(\frac{1}{q_m^{2s}(\bfa) B^{mg(s)}}\right) \leq 1
\right\}.
\end{equation}
As in \cite[Section 2]{WangWu2008}, we may show that
\begin{equation}\label{Eq:LimSmB}
\lim_{m\to\infty} s_{m,B}
=
s_B.
\end{equation}

The lower estimate for the Hausdorff dimension of sets resembling that in Theorem \ref{TEO:MAIN:HD} tends to be computationally demanding. 
A general strategy is to construct a Cantor set contained in the set of interest, define a probability measure on it, and apply the Mass Distribution Principle \cite[Proposition 4.2]{FalconerBookFractal}. 
Recently, Hussain and Shulga \cite{ShulgaHussain2023Exp} computed the Hausdorff dimension of a large family of sets defined by continued fractions restrictions. 
For $\ell\in\Na$ and $A_0,\ldots, A_{\ell-1}>1$, write  
\begin{equation}\label{Eq:Def:SetS}
S_{\ell}(A_0,\ldots, A_{\ell-1})
:=
\left\{ 
x\in [0,1)\setminus\QU: 
\begin{matrix}
A_0^n\leq a_n(x)\leq 2A_0^n, \ldots, A_{\ell-1}^n\leq a_{n+\ell-1} (x)\leq 2A_{\ell-1}^n \\ \text{ for infinitely many } n\in\Na
\end{matrix}
\right\}.
\end{equation}

Put $\beta_{-1}=1$ and $\beta_{i}=A_i\beta_{i-1}$ for $i\in \{0,\ldots, \ell-1\}$ and, for such $i$, put
\[
d_i
:=
\inf\left\{ s\geq 0: P(T,-s\log|T'| - s\log \beta_i + (1-s)\log\beta_{i-1} )\leq 0\right\}.
\]
\begin{lem01}[{\cite[Theorem 1.1]{ShulgaHussain2023Exp}}]\label{Le:MumtazNikita}
$\hdim S_{\ell}(A_0,\ldots, A_{\ell-1})
=
\displaystyle\min_{0\leq i\leq \ell-1} d_i$.
\end{lem01}
\section{Series estimates}\label{Section:Series}
\subsection{Bounds for the Lebesgue measure}
In this section, we estimate the series needed in the proofs of theorems \ref{TEO:MAIN:MEASURE} and \ref{TEO:MAIN:HD}.
The variables $a_1,a_2,\ldots$, $b_1, b_2, \ldots$, $c_1, c_2, \ldots$ take values in the natural numbers.

We start with the lemmas needed in the proof of Theorem \ref{TEO:MAIN:MEASURE}. 
The first can be retrieved from the proof of \cite[Theorem 1.5]{HuangWuXu2020} and the second is a well-known result.
\begin{lem01}[\cite{HuangWuXu2020}]\label{Lem:HWX}
Given $j\in\Na$,  for $M>1$ we have
\begin{equation}\label{Eq:Series:01}
\sum_{a_1\cdots a_{j}\geq M }
\frac{1}{a_1^2\cdots a_{j}^2} 
\asymp_{j}
\frac{\log^{j-1}(M)}{M}.
\end{equation}
\end{lem01}


\begin{lem01}[Dirichlet-Piltz formula]
Take $k\in\Na$. For each $x\in \Na$, let $d_k(x)$ denote the number of ways we can express $x$ as a product of $k$ natural numbers. For every $M>1$, define
\[
D_k(M):= \sum_{x\leq M} d_k(x).
\]
There exist a polynomial $P_k$ of degree $k-1$ and a constant $0<\alpha_k<1$ such that for any $\veps>0$ we have
\[
D_k(M)
=
MP_k(\log M) + O(M^{\alpha_k+\veps})
\;\text{ as }\; M\to\infty.
\]
\end{lem01}
\begin{proof}
See, for example, \cite[Section 4]{HardyLittlewood1923}
\end{proof}

Clearly, for any $\alpha>0$ we have
\begin{align}
\int \frac{\log^n(\alpha t)}{t} \md t =  \frac{\log^{n+1}(\alpha t)}{n+1} + C \quad \text{for all } n\in\Na, \label{Eq:IntLog:01} \\
\int \log(\alpha t) \md t =  t\log(\alpha t) - t +C . \label{Eq:IntLog:02}
\end{align}

\begin{lem01}\label{Lem:SerEstNu}
For $\ell\in\Na$ and $M>1$, we have
\begin{equation}\label{Eq:Series:03}
\sum_{a_1\cdots a_{\ell}\leq M} \frac{1}{a_1\cdots a_{\ell}} 
\asymp_{\ell}
\log^{\ell}(M), 
\end{equation}
\begin{equation}  \label{Eq:Series:04}
\sum_{\substack{a_1\cdots a_{\ell}<M \\ a_2\cdots a_{\ell+1}\geq M}} \frac{1}{a_1^2\cdots a_{\ell}^2a_{\ell+1}^2 }
\asymp_{\ell}
\frac{\log^{\ell-1}M}{M}. 
\end{equation}
\end{lem01}
\begin{proof}
The estimate \eqref{Eq:Series:03} is clear when $\ell=1$. 
If it holds for some $\ell\in\Na$, then, by \eqref{Eq:IntLog:01},
\[
\sum_{a_1\cdots a_{\ell}a_{\ell+1}\leq M} \frac{1}{a_1\cdots a_{\ell}a_{\ell+1}} 
=
\sum_{a_{\ell+1}=1 }^M
\frac{1}{a_{\ell+1}} 
\sum_{a_1\cdots a_{\ell} \leq \frac{M}{a_{\ell+1}}} 
\frac{1}{a_1\cdots a_{\ell}} 
\asymp
\sum_{a_{\ell+1}=1 }^M
\frac{1}{a_{\ell+1}} 
\log^{\ell}\left( \frac{M}{a_{\ell+1}}\right)
\asymp
\log^{\ell+1} M.
\]
For \eqref{Eq:Series:04}, note that
\[
\sum_{\substack{a_1\cdots a_{\ell}<M \\ a_2\cdots a_{\ell+1}\geq M}} \frac{1}{a_1^2\cdots  a_{\ell+1}^2 }
=
\sum_{a_2\cdots a_{\ell}\leq M}
\sum_{ a_1\leq \frac{M}{a_2\ldots a_{\ell}} } 
\sum_{a_{\ell+1} \geq \frac{M}{a_2\cdots a_{\ell}}} 
\frac{1}{a_1^2  \cdots a_{\ell+1}^2}.
\] 
We split the series to the right by conditioning on whether $1\leq a_2\ldots a_{\ell}\leq \frac{M}{2}$ or $\frac{M}{2}< a_2\cdots a_{\ell} < M$. 
In the first case, since super-harmonic series are convergent, we have
\begin{align*}
\sum_{a_2\cdots a_{\ell}\leq \frac{M}{2} }
\sum_{ a_1\leq \frac{M}{a_2\ldots a_{\ell}} } 
\sum_{a_{\ell+1} \geq \frac{M}{a_2\cdots a_{\ell}}} 
\frac{1}{a_1^2  \cdots a_{\ell+1}^2}
&= 
\sum_{a_2\cdots a_{\ell}\leq \frac{M}{2} } 
\frac{1}{a_2^2\cdots a_{\ell}^2}
\sum_{ a_1\leq \frac{M}{a_2\ldots a_{\ell+1}} } 
\frac{1}{a_1^2}
\sum_{a_{\ell+1} \geq \frac{M}{a_2\cdots a_{\ell}}} 
\frac{1}{a_{\ell+1}^2} \\
&\asymp 
\sum_{a_2\cdots a_{\ell}\leq \frac{M}{2} }
\frac{1}{a_2^2 \cdots a_{\ell}^2}
\left( \sum_{a_{\ell+1} \geq \frac{M}{a_2\cdots a_{\ell}}} \frac{1}{a_{\ell+ 1}^2}\right)   \\
&\asymp
\frac{1}{M}
\sum_{a_2\cdots a_{\ell}\leq \frac{M}{2} }
\frac{1}{a_2 \cdots a_{\ell}} \\
&\asymp
\frac{\log^{\ell-1}M}{M} &&\text{(by \eqref{Eq:Series:03})}.
\end{align*}
When $\frac{M}{2}< a_2\cdots a_{\ell} < M$, the condition $a_1\leq \frac{M}{a_2\ldots a_{\ell}}$ reduces to $a_1=1$ and $a_{\ell+1} \geq \frac{M}{a_2\cdots a_{\ell}}$ becomes $a_{\ell+1}\geq 2$, then
\begin{align*}
\sum_{\frac{M}{2}\leq a_2\cdots a_{\ell}\leq M }
\sum_{ a_1\leq \frac{M}{a_2\ldots a_{\ell}} } 
\sum_{a_{\ell+1} \geq \frac{M}{a_1\cdots a_{\ell}}} 
\frac{1}{a_1^2 \cdots a_{\ell+1}^2}
&\asymp
\sum_{\frac{M}{2}\leq a_2\cdots a_{\ell}\leq M }
\frac{1}{a_2^2 \cdots a_{\ell}^2} \\
&\ll
\sum_{\frac{M}{2}\leq a_2\cdots a_{\ell} }
\frac{1}{a_2^2 \cdots a_{\ell}^2}  \\
&\ll
\frac{\log^{\ell-2}M}{M}. &&\text{(by Lemma \ref{Lem:HWX})}.
\end{align*}

This proves \eqref{Eq:Series:04}.
\end{proof}
Theorem \ref{Thm:OverlapEstimatesSeries} below is the main technical tool in the proof of Theorem \ref{TEO:MAIN:MEASURE}.
We will use it to estimate the measure of the sets where the restricted blocks of consecutive partial quotients overlap.
\begin{teo01}\label{Thm:OverlapEstimatesSeries}
For any $r,j\in\Na$, we have
\[
\sum_{\substack{a_1 \cdots a_r b_1\cdots b_j\geq M \\ b_1\cdots b_j c_1 \cdots c_r\geq M}} 
\frac{1}{a_1^2\cdots a_r^2 b_1^2\cdots b_j^2 c_1^2 \cdots c_r^2}
\asymp_{r}
\frac{\log^{j-1}(M)}{M}
\quad\text{ as }M\to\infty.
\]
\end{teo01}
The following lemma contains the core of the proof of Theorem \ref{Thm:OverlapEstimatesSeries}.

\begin{lem01}\label{Lem:OverlapEstimatesSeries:bis}
For each $r,j\in\Na$, there are polynomials $P_{r,j}, Q_{r,j}\in\RE[X]$ with $Q_{r,1}=0$ and $\deg Q_{r,j} = j-2$ if $j\geq 2$, and a positive constant $c_{r,j}$ such that
\[
\int_{\substack{1\leq x_1\cdots x_j\leq M\\ 1\leq x_1, \ldots, 1\leq x_j}}
\log^{2(r-1)}\left(\frac{x_1\cdots x_j}{M}\right)
\md \bfx
=
c_{r,j} M\log^{j-1} M +  P_{r,j}(\log M) + MQ_{r,j}(\log M).
\]
\end{lem01}
An immediate consequence of Lemma \ref{Lem:OverlapEstimatesSeries:bis} is the next corollary.

\begin{coro01}\label{Coro:OverlapEstimatesSeries}
For each $r,j\in\Na$, 
\[
\frac{1}{M^2}
\int_{\substack{1\leq x_1\cdots x_j\leq M\\ 1\leq x_1, \ldots 1\leq x_j}}
\log^{2(r-1)}\left(\frac{x_1\cdots x_j}{M}\right)
\md x
\asymp_{r,j}
\frac{\log^{j-1}(M)}{M}
\quad\text{ as }M\to\infty.
\]
\end{coro01}
Lemmas \ref{Lem:r=1} and \ref{Lem:j=1} below are intermediate steps toward the proof of Lemma \ref{Lem:OverlapEstimatesSeries:bis}.
\begin{lem01}\label{Lem:r=1}
For each $j\in\Na$, there is a polynomial $R_j\in \RE[X]$ satisfying $R_1=0$ and $\deg R_j=j-2$ for $j\geq 2$ such that
\[
\int_{\substack{1\leq x_1\cdots x_j\leq M \\ 1\leq x_1, \ldots, 1\leq x_j}} 1 \,\md \bfx
=
\frac{1}{(j-1)!} M\log^{j-1}M + MR_j(\log M) + (-1)^j.
\]
\end{lem01}
\begin{proof}
The case $j=1$ is trivial.
For $j=2$, we have
\[
\int_{\substack{1\leq x_1 x_2 \leq M\\ 1\leq x_1, 1\leq x_{2}}}
1
\md \bfx
=
\int_1^M
\left(
\int_{1}^{\frac{M}{x_2}}
1
\md x_1
\right)
\md x_{2}  
=
\int_1^M
\frac{M}{x_2} - 1
\md x_{2}  
=
M\log M  - M + 1,
\]
so $R_2(X)=-1$. 
Assume the result holds for some $j\in\Na$; then,
\begin{align*}
\int_{\substack{1\leq x_1\cdots x_jx_{j+1}\leq M\\ 1\leq x_1, \ldots, 1\leq x_{j+1}}}
1
\md \bfx
&=
\int_1^M
\left(
\int_{\substack{1\leq x_1\cdots x_j\leq \frac{M}{x_{j+1}} \\ 1\leq x_1, \ldots, 1\leq x_j}}
1
\md(x_1,\ldots, x_j)
\right)
\md x_{j+1} \nonumber\\
&=
\int_1^M
\frac{1}{(j-1)!} \frac{M}{x_{j+1}} \log^{j-1}\left(\frac{M}{x_{j+1}}\right)
+
\frac{M}{x_{j+1}} R_{j}\left( \log \frac{M}{x_{j+1}} \right) 
+
(-1)^j
\md x_{j+1} \nonumber\\
&=
\frac{1}{j!} M\log^{j} M
+
MR_{j+1}(\log M)
+
(-1)^{j+1} \quad\quad \text{(by \eqref{Eq:IntLog:01} and \eqref{Eq:IntLog:02}).}
\end{align*}
Finally, note that $\deg R_{j+1} = 1+ \deg R_j = (j+1)-2$.
\end{proof}
\begin{lem01}\label{Lem:j=1}
For $k\in\Na$, there is a polynomial $\tilde{R}_k\in \RE[X]$ of degree $\deg \tilde{R}_k=k$ and leading coefficient $(-1)^{k+1}$ such that for all $M>1$ we have
\[
\int_1^M \log^k\left( \frac{x_1}{M}\right) \md x_1
=
(-1)^k k! M+\tilde{R}_k(\log M).
\]
In particular, for each $r\in \Na$ there is a polynomial $\tilde{R}_{2(r-1)}\in\RE[x]$ of degree $2(r-1)$ and leading coefficient $-1$ such that
\[
\int_1^M \log^{2(r-1)}\left( \frac{x_1}{M}\right) \md x_1
=
(2r-2)! M+\tilde{R}_{2(r-1)}(\log M).
\]
\end{lem01}
\begin{proof}
For the base $k=1$, we use  \eqref{Eq:IntLog:02} to obtain
\[
\int_1^M \log\left( \frac{x_1}{M}\right) \md x_1
=
-M + \log M +1.
\]
Note that $\tilde{R}_1(X)=X$.
Assume the result holds for some $k\in\Na$.
Integrating by parts, 
\begin{align*}
\int_1^M \log^{k+1}\left( \frac{x_1}{M}\right) \md x_1
&=
(-1)^{k+2}\log^{k+1}(M)
-
(k+1)\int_1^M \log^k\left(\frac{x_1}{M}\right) \md x_1 \\
&=
(-1)^{k+2}\log^{k+1}(M)
+
(-1)^{k+1}(k+1)!
M
-
(k+1)\tilde{R}_k(\log M) \\
&=
(-1)^{k+1}(k+1)!
M
+
\tilde{R}_{k+1}(\log M)
\end{align*}
and $\tilde{R}_{k+1}(X)=(-1)^{k+2}X^{k+1} - (k+1)\tilde{R}_k(X)$.
\end{proof}
\begin{proof}[Proof of Lemma \ref{Lem:OverlapEstimatesSeries:bis}]
We show the lemma by induction on $r$. 
The base $r=1$ and $j\in\Na$ is precisely Lemma \ref{Lem:r=1}.
Assuming that the result holds for $r\in\Na$ and all $j\in\Na$, we prove it for $r+1$ and all $j$.
The case $j=1$ follows directly from Lemma \ref{Lem:j=1}.
For $j=2$, Lemma \ref{Lem:j=1} yields
\begin{align*}
\int_{\substack{1\leq x_1 x_{2}\leq M \\ 1\leq x_1, 1\leq x_{2}}}
\log^{2r}\left( \frac{x_1x_2 }{M}\right) \md \bfx
&= 
\int_1^M
\int_{1}^{\frac{M}{x_2}} 
\log^{2r}\left( \frac{x_1 }{\frac{M}{x_{2}}}\right) 
\md x_1
\md x_{2}
\\
&=
\int_1^M
(2r)! \frac{M}{x_2} +\tilde{R}_{2r}\left(\log \frac{M}{x_2}\right)
\md x_{2} \\
&=
(2r)! M
\int_1^M
\frac{1}{x_2}\md x_2 
+
\int_1^M
\tilde{R}_{2r}\left(\log \frac{M}{x_2}\right)
\md x_{2} \\
&=
(2r)! M\log M
+
\int_1^M
\tilde{R}_{2r}\left(\log \frac{M}{x_2}\right)
\md x_{2}
\end{align*}
We apply Lemma \ref{Lem:j=1} on each term of $\tilde{R}_{2r}\left(\log \frac{M}{x_2}\right)$ to obtain polynomials $Q_{r+1,2}, P_{r+1, 2}\in \RE[X]$ such that $Q_{r+1,2}$ is constant and
\[
\int_{\substack{1\leq x_1 x_{2}\leq M \\ 1\leq x_1, 1\leq x_{2}}}
\log^{2r}\left( \frac{x_1x_2 }{M}\right) \md \bfx
=
(2r)! M\log M + P_{r+1,2}(\log M) + MQ_{r+1,2}(\log M).
\]
If the result holds for $j\in\Na$, then
\begin{align*}
\int_{\substack{1\leq x_1\cdots x_{j+1}\leq M \\ 1\leq x_1, \ldots, 1\leq x_{j+1}}}
&\log^{2r}\left( \frac{x_1\cdots x_jx_{j+1}}{M}\right) \md \bfx
= \\
&=
\int_1^M
\left(
\int_{\substack{1\leq x_1\cdots x_{j}\leq \frac{M}{x_{j+1}} \\ 1\leq x_1, \ldots, 1\leq x_{j}}}
\log^{2r}\left( \frac{x_1\cdots x_j }{\frac{M}{x_{j+1}}}\right) 
\md (x_1, \ldots, x_j)
\right)
\md x_{j+1} \\
&=
\int_1^M
c_{r,j}\frac{M}{x_{j+1}}\log^{j-1}\left( \frac{M}{x_{j+1}}\right)
+
P_{r+1, j}\left(\log  \frac{M}{x_{j+1}} \right)
+
\frac{M}{x_{j+1}} Q_{r+1, j}\left(\log \frac{M}{x_{j+1}} \right)
\md x_{j+1} \\
&=
\frac{c_{r,j}}{j} M \log^{j}M
+
P_{r+1, j+1}(\log M)
+
MQ_{r+1, j+1}(\log M).
\end{align*}
By Lemma \ref{Lem:r=1}, $\deg P_{r+1, j+1} = 1 + \deg P_{r+1, j}$ and, by \eqref{Eq:IntLog:02}, the polynomial $Q_{r+1, j+1}$ satisfies $\deg Q_{r+1, j+1} = 1 + \deg Q_{r+1, j}$.
Moreover, in view of Lemma \ref{Lem:r=1}, when computing the integral related to $P_{r+1, j}$, a scalar multiple of $M$ appears, and we let $M Q_{r+1, j+1}(\log M)$ absorb it.
This proves the lemma.
\end{proof}
\begin{proof}[Proof of Theorem \ref{Thm:OverlapEstimatesSeries}]
First, assume that $r=1$ and note that 
\[
\sum_{\substack{a_1 b_1\cdots b_j\geq M \\ b_1\cdots b_j c_1 \geq M } }
\frac{1}{a_1^2 b_1^2\cdots b_j^2 c_1^2 }
=
\sum_{\substack{a_1  b_1\cdots b_j\geq M \\ b_1\cdots b_j c_1  \geq M\\ b_1\cdots b_j\geq M} }
\frac{1}{a_1^2 b_1^2\cdots b_j^2 c_1^2 }
+
\sum_{\substack{a_1 b_1\cdots b_j\geq M \\ b_1\cdots b_j c_1 \geq M\\ b_1\cdots b_j< M}} 
\frac{1}{a_1^2 b_1^2\cdots b_j^2 c_1^2 }.
\]
The first sum to the right is equivalent, up to multiplicative constants, to $M^{-1} \log^{j-1}(M)$ by Lemma \ref{Lem:HWX}, because there are no actual restrictions on the natural numbers $a_1$ and $c_1$. 
For the second term, we have
\[
\sum_{b_1\cdots b_j< M} 
\frac{1}{b_1^2\cdots b_j^2}
\left( \sum_{a> \frac{M}{b_1\cdots b_j}} \frac{1}{a^2}\right)^2
\ll 
\sum_{b_1\cdots b_j< M} 
\frac{1}{b_1^2\cdots b_j^2}
\frac{1}{\frac{M^2}{b_1^2\cdots b_j^2}}
=
\frac{1}{M^2}
\sum_{b_1\cdots b_j< M} 
1
\asymp
\frac{\log^{j-1} M}{M},
\]
where the last estimate follows from the Dirichlet-Piltz formula.
Now, let's assume that $r\geq 2$.
We split the series into two sums according to the size of $b_1\cdots b_j$.
When $b_1\cdots b_j\geq M$, we have
\[
\sum_{\substack{a_1 \cdots a_r b_1\cdots b_j\geq M \\ b_1\cdots b_j c_1 \cdots c_r\geq M\\ b_1\cdots b_j\geq M}} 
\frac{1}{a_1^2\cdots a_r^2 b_1^2\cdots b_j^2 c_1^2 \cdots c_r^2}
=
\sum_{b_1\cdots b_j\geq M}
\frac{1}{b_1^2\cdots b_j^2} \left(\sum_{a\in\Na}\frac{1}{a^2} \right)^{2r}
\asymp_r
\frac{\log^{j-1}M}{M}.
\]

For the case $b_1\cdots b_j< M$, we further divide the corresponding sum by considering whether the additional restriction $b_1\cdots b_j=1$ is true or false:
\[
\sum_{\substack{a_1 \cdots a_r b_1\cdots b_j\geq M \\ b_1\cdots b_j c_1 \cdots c_r\geq M\\ b_1\cdots b_j =1 }} 
\frac{1}{a_1^2\cdots a_r^2 b_1^2\cdots b_j^2 c_1^2 \cdots c_r^2}
=
\left(
\sum_{\substack{a_1 \cdots a_r \geq M }} 
\frac{1}{a_1^2\cdots a_r^2}
\right)^2
\asymp
\frac{\log^{2(r-1)}M}{M^2}.
\]
If $b_1\cdots b_j>1$, then $b_k\geq 2$ for some $k$, and
\begin{align}
\sum_{\substack{a_1 \cdots a_r b_1\cdots b_j\geq M \\ b_1\cdots b_j c_1 \cdots c_r\geq M\\ 2\leq b_1\cdots b_j < M}} 
\frac{1}{a_1^2\cdots a_r^2 b_1^2\cdots b_j^2 c_1^2 \cdots c_r^2} 
&=
\sum_{2\leq  b_1\cdots b_j<M} \frac{1}{b_1^2\cdots b_j^2} \left( \sum_{a_1\cdots a_r\geq \frac{M}{b_1\cdots b_j}} \frac{1}{a_1^1\cdots a_r^2}\right)^2 \nonumber \\
&\asymp
\frac{1}{M^2}
\sum_{2\leq  b_1\cdots b_j<M} 
\log^{2(r-1)} \left( \frac{M}{b_1\cdots b_j}\right). \label{Eq:BoundIntegral}
\end{align}
Let $\clE\subseteq \Na^j$ be the set of vectors $\bfb=(b_1, \ldots, b_j)\in\Na^{j}$ such that $2\leq b_1\cdots b_j\leq M$.
For $\bfb\in \clE$, define
\[
V(\bfb)
=
\prod_{k=1}^j [b_k, b_k+1].
\]
Then, for all $\bfx=(x_1,\ldots, x_j)\in V(\bfb)$, we have $b_1\cdots b_j\leq x_1\cdots x_j\leq 2^j b_1\cdots b_j$, so
\[
\log\left(\frac{b_1\cdots b_j}{2^j M}\right) 
\leq 
\log\left(\frac{x_1\cdots x_j}{2^j M}\right) 
\leq 
\log\left(\frac{b_1\cdots b_j}{M}\right) 
\leq 
0,
\]
and, since $2(r-1)$ is a positive even integer,
\[
0\leq 
\log^{2(r-1)}\left(\frac{b_1\cdots b_j}{M}\right) 
\leq 
\log^{2(r-1)}\left(\frac{x_1\cdots x_j}{2^j M}\right) 
\leq 
\log^{2(r-1)}\left(\frac{b_1\cdots b_j}{2^j M}\right) .
\]
Furthermore, if $E=\bigcup_{\bfb\in \clE} V(\bfb)$, we have
\[
E
\subseteq 
F:=
\left\{ (x_1, \ldots, x_j)\in \RE^j: 1\leq x_1, \ldots, 1\leq x_j \text{ and } 1\leq x_1\cdots x_j\leq 2^jM\right\};
\]
hence, writing $\bfx=(x_1, \ldots, x_j)$,
\begin{align*}
\frac{1}{M^2}
\sum_{\bfb\in\clE} \log^{2(r-1)}\left( \frac{b_1\cdots b_j}{M}\right) 
&\leq 
\frac{1}{M^2}
\int_E \log^{2(r-1)}\left( \frac{x_1\cdots x_j}{2^jM}\right) \md \bfx \\ 
&\leq 
\frac{1}{M^2}
\int_{\substack{1\leq x_1, \ldots, 1\leq x_j \\ 1\leq x_1 \cdots x_j \leq 2^j M }} \log^{2(r-1)}\left( \frac{x_1\cdots x_j}{2^jM}\right) \md \bfx \\ 
&\asymp_r \frac{\log^{j-1}(2^jM)}{ M} &&{\text{(by Corollary \ref{Coro:OverlapEstimatesSeries})}} \\
&\asymp_{r,j} \frac{\log^{j-1}(M)}{M}.
\end{align*} 
\end{proof}
\subsection{Bounds for the Hausdorff dimension}
Theorem \ref{TEO:MAIN:HD} requires similar yet simpler estimates.

 \begin{lem01}[{\cite[Lemma 4.2]{HuangWuXu2020}}]\label{Lem:HWX:s:estimate}
For $\beta>1$, $0<s<1$, and $n\in\Na$, we have
\begin{equation}\label{Eq:Series:05}
\sum_{ a_1\cdots a_{\ell} \leq \beta^n}
\left(\frac{1}{a_1\cdots a_{\ell}} \right)^s
\asymp
\frac{(\log \beta^n)^{{\ell}-1}}{({\ell}-1)!} \beta^{n(1-s)}.  
\end{equation}
\end{lem01}
\begin{lem01}\label{Lem:s:estimate:02}
If $M, t>1$, then
\begin{equation}\label{Eq:Series:06}
\sum_{ab\geq M} \frac{1}{a^tb^t}
\asymp _t
\frac{1}{t-1}\,\frac{\log M}{M^{t-1}},
\end{equation}
\begin{equation}\label{Eq:Series:07}
\sum_{abc\geq M} \frac{1}{a^tb^tc^t}
\asymp_t
\frac{\log^2 M}{M^{t-1}}.
\end{equation}
\end{lem01}
\begin{proof}
To show \eqref{Eq:Series:06}, we separate the series into two sums according to whether $a\geq M$ or $1\leq a<M$ holds.
In the first case, we obtain a series of order $\frac{1}{(t-1)M^{t-1}}$.
The sum corresponding to $1\leq a<M$ has order $\frac{\log M}{(t-1)M^{t-1}}$.
We use a similar strategy on \eqref{Eq:Series:07}. 
First, we split the sum as follows:
\[
\sum_{abc\geq M} \frac{1}{a^tb^tc^t}
=
\sum_{\substack{abc\geq M \\  c\geq M} } \frac{1}{a^tb^tc^t}+
\sum_{\substack{abc\geq M \\ c< M} } \frac{1}{a^tb^tc^t}.
\]
When $c\geq M$, the parameters $a,b$ can take any value in $\Na$, so
\[
\sum_{\substack{abc\geq M \\  c\geq M} } \frac{1}{a^tb^tc^t}
\asymp 
\frac{1}{(t-1)^3} \frac{1}{M^{t-1}}.
\]
When $c < M$, we have 
\begin{align*}
\sum_{\substack{abc\geq M \\ c< M} } 
\frac{1}{a^tb^tc^t}
&=
\sum_{1\leq c< M} 
\frac{1}{c^t}
\sum_{ ab\geq \frac{M}{c}}
\frac{1}{a^tb^t} \\
&\asymp 
\frac{1}{t-1}
\sum_{1\leq c< M} 
\frac{1}{c^t}
\left(\frac{c}{M}\right)^{t-1}
\log \frac{M}{c} \\
&=
\frac{1}{(t-1)M^{t-1}}
\sum_{1\leq c< M} 
\frac{1}{c}
\log \frac{M}{c} \\
&\asymp
\frac{1}{(t-1)}
\frac{\log^2 M}{M^{t-1}} &&\text{(by \eqref{Eq:IntLog:02})}.
\end{align*}
\end{proof}

\begin{lem01}\label{Lem:EstseriesRealExp}
For any $\frac{1}{2}<s<1$ and $M>1$ we have
\[
\sum_{ab\leq M} \frac{1}{a^sb^s} \asymp_s M^{1-s} \log M.
\]
\end{lem01}
\begin{proof}
For all $M>1$, we have
\begin{equation}\label{Lem:EstSerieSup:01}
\int_{\substack{1\leq x, 1\leq y \\ 1\leq xy \leq M}} \frac{1}{x^sy^s} \md x \md y 
=
\frac{M^{1-s}}{1-s} \log M + \frac{1}{(1-s)^2} \left( 1 - M^{1-s}\right).    
\end{equation}
For every $N>1$ define the set 
\[
H_N=\left\{(x, y)\in\RE^2: 1\leq x, \, 1\leq y, \, xy\leq M\right\}.
\]
Let us check that 
\[ 
\int_{H_M} \frac{1}{x^sy^s} \md x \md y
\asymp_s
\sum_{ab\leq M} \frac{1}{a^s b^s }
\]
Define $L_M=H_M\cap \Za^2$ and, for each $(a,b)\in L_M$, put $V(a,b)=[a, a+1)\times [b,b+1)$.
Let $f,g:\RE_{\geq 1} \times \RE_{\geq 1}\to\RE$ be given by 
\[
f(x,y) = \frac{1}{x^sy^s}, 
\quad 
g(x,y) = \frac{1}{[x]^s[y]^s},
\]
so $f(x,y)\leq g(x,y)$ for all $(x,y)$. 
Also, for any $(a,b)\in L_M$, we have 
\[
\int_{V(a,b)} f(x,y)\md x \md y 
\leq 
\int_{V(a,b)} g(x,y)\md x \md y 
=
\frac{1}{a^sb^s}.
\]
Therefore, writing $\cL_M = \bigcup_{(a,b)\in L_M} V(a,b)$, we have $H_M\subseteq \clL_M$ and
\[
\int_{H_M} f(x,y) \md x \md y 
\leq 
\int_{\cL_M} f(x,y)\md x \md y 
\leq 
\int_{\cL_M} g(x,y)\md x \md y 
=
\sum_{ab\leq M} \frac{1}{a^sb^s}.
\]
For the other inequality, consider the sets
\[
A = \{(a,b)\in L_M: \min\{a,b\}\geq 2\}, \quad
\text{ and }
B= L_M \setminus A.
\]
On the one hand, $\displaystyle\sum_{(a,b)\in B} a^{-s}b^{-s} \asymp M^{1-s}$. 
On the other hand, for all $(a,b)\in A$ we have $V(a-1,b-1)\subseteq H_M$, so 
\begin{align*}
\sum_{(a,b)\in A} \frac{1}{a^s b^s}
&=
\sum_{(a,b)\in A}
\int_{V(a,b) } g(x,y)\md x \md y \\ 
&=
\sum_{(a,b)\in A}
\int_{V(a-1,b-1) } f(x,y)\md x \md y  
\leq 
\int_{ H_M } f(x,y)\md x \md y.
\end{align*}
\end{proof}

\section{Proof of Theorem \ref{TEO:MAIN:MEASURE}}\label{Section:ProofMeasure}
For $\ell \in\Na_{\geq 3}$, $n\in\Na_{\geq \ell}$, and $k\in \{1,\ldots, n-1\}$, define
\[
A_{n,k}(\vphi)
:=
\left\{x\in [0,1) : 
\begin{matrix}
a_k(x)\cdots a_{k+\ell - 1}(x) \geq \vphi(n), \\
a_n(x)\cdots a_{n+\ell - 1}(x) \geq \vphi(n)
\end{matrix}
 \right\}
\]
and write
\[
A_n(\vphi)
:=
\bigcup_{k=1}^{n-1} A_{n,k}(\vphi).
\]
As in \cite{TanZhou2024}, the monotonicity of $\vphi$ implies
\[
\clF_{\ell}(\vphi)
=
\limsup_{n\to\infty} A_n(\vphi).
\]
\subsection{Convergence case}
We show that 
\begin{equation}\label{Eq:ConverGen}
\sum_{n\geq 1}
\left(\frac{n\log^{2(\ell-1)}\vphi(n)}{\vphi^2(n)} + \frac{\log^{\ell-2}\vphi(n)}{\vphi(n)}\right)<\infty
\quad\text{ implies }\quad
\leb(\clF_{\ell}(\vphi))=0.
\end{equation}
We estimate $\leb(A_{n,k})$ by considering two cases according to whether the restricted blocks overlap or not.
The contribution of non-overlapping blocks with large products is of order
\[
  \frac{n\log^{2(\ell-1)}\vphi(n)}{\vphi^2(n)}.
\]
The contribution of overlapping blocks with large products is of order
\[
  \frac{\log^{\ell-2}\vphi(n)}{\vphi(n)}.
\]

\begin{lem01}\label{Lem:Leb:Convergence}
Let $k\in\{1,\ldots, n-1\}$.
\begin{enumerate}[1.]
\item If $1\leq k \leq n-\ell$, then 
\[
\leb( A_{n,k}(\vphi))
\asymp_{\ell} 
\frac{\log^{2(\ell-1)}\vphi(n)}{\vphi^2(n)}.
\]
\item If $n-\ell +1 \leq k\leq n- 1$, then 
\[
\leb( A_{n,k}(\vphi))
\asymp_{\ell} 
\frac{\log^{ k+\ell-n-1 }\vphi(n)}{\vphi(n)}.
\]
\end{enumerate}
\end{lem01}
\begin{proof}
\begin{enumerate}[1.]
\item Assume that $1\leq k\leq n-\ell$, so $k+\ell-1 \leq n-1$. Let $\bfb=(b_1,\ldots, b_{n+\ell-1})\in\Na^{n+\ell-1}$ be such that $A_{n,k}(\vphi)\cap I(\bfb)\neq \vac$, so $I(\bfb)\subseteq A_{n,k}(\vphi)$. 
Write
\[
\bfb'=(b_1,\ldots, b_{k-1},b_{k+\ell},\ldots, b_{n-1});
\]
that is, $\bfb'$ is $\bfb$ without the blocks $(b_{k},\ldots, b_{k+\ell-1})$ and $(b_n,\ldots, b_{n+\ell-1})$, then
\[
\leb(I_{n+\ell-1}(\bfb))
\asymp_{\ell}
\leb(I_{n-\ell-1}(\bfb')) \, 
\frac{1}{(b_{k}^2\cdots b_{k+\ell-1}^2) (b_{n}^2\cdots b_{n+\ell-1}^2)}.
\]
Therefore, we have
\begin{align*}
\leb(A_{n,k}(\vphi))
&=
\sum_{\bfb\in \Na^{n+\ell-1}} \leb\left(A_{n,k}(\vphi)\cap I_{n+\ell-1}(\bfb)\right)
\nonumber\\
&=
\sum_{\substack{\bfb\in \Na^{n+\ell-1}\\ b_{k}\cdots b_{k+\ell-1}\geq \vphi(n)\\ b_{n}\cdots b_{n+\ell-1}\geq \vphi(n) }} 
\leb\left(I_{n+\ell-1}(\bfb)\right) \nonumber\\
&\asymp_{\ell}
\sum_{\substack{\bfb'\in \Na^{n-\ell-1}\\ b_{k}\cdots b_{k+\ell-1}\geq \vphi(n)\\ b_{n}\cdots b_{n+\ell-1}\geq \vphi(n) }} 
\leb(I_{n-\ell-1}(\bfb')) \, 
\frac{1}{(b_{k}^2\cdots b_{k+\ell-1}^2) (b_{n}^2\cdots b_{n+\ell-1}^2)}
\nonumber\\
&=
\left(\sum_{\bfb'\in \Na^{n-\ell-1}} 
\leb(I_{n-\ell-1}(\bfb') ) \right)
\left(\sum_{b_{k}\cdots b_{k+\ell-1}\geq \vphi(n) } 
\frac{1}{b_{k}^2\cdots b_{k+\ell-1}^2}\right)
\left(\sum_{b_{n}\cdots b_{n+\ell-1}\geq \vphi(n) } 
\frac{1}{b_{n}^2\cdots b_{n+\ell-1}^2}\right) \nonumber\\
&=
\left(\sum_{b_{1}\cdots b_{\ell}\geq \vphi(n) } 
\frac{1}{b_{1}^2\cdots b_{\ell}^2}\right)^2  \\
&\asymp_{\ell}
\frac{\log^{2(\ell-1)}\vphi(n)}{\vphi^2(n)}.
\end{align*}
The last estimate follows from Lemma \ref{Lem:HWX}.
\item Assume that $n-\ell +1 \leq k\leq n- 1$ and write $j=k+ \ell -n$ and $r= \ell-j$.
Then, $j$ is the length of the overlap of the blocks mentioned in the definition of $A_{n,k}(\vphi)$ and $r$ is the length of the non-overlapping blocks. 
Writing $\bfb'= (b_1, \ldots, b_{k-1})$ and noting that $n+j+r-1=n+\ell-1$ and $k+r-1=n-1$, we have
\begin{align*}
\leb(A_{n,k}(\vphi))
&=
\sum_{\bfb\in \Na^{n+\ell-1}} \leb\left(A_{n,k}(\vphi)\cap I_{n+\ell-1}(\bfb)\right)
\nonumber\\
&=
\sum_{\substack{\bfb\in \Na^{n+\ell-1}\\ b_{k}\cdots b_{k+\ell-1}\geq \vphi(n)\\ b_{n}\cdots b_{n+\ell-1}\geq \vphi(n) }} 
\leb\left(I_{n+\ell-1}(\bfb)\right) \nonumber\\
&\asymp_{\ell}
\sum_{\substack{\bfb'\in \Na^{n-\ell-1}\\ b_{k}\cdots b_{k+\ell-1}\geq \vphi(n)\\ b_{n}\cdots b_{n+\ell-1}\geq \vphi(n) }} 
\leb(I_{k-1}(\bfb')) \, 
\frac{1}{ b_{k}^2\cdots b_{k+r-1}^2  b_{n}^2 \cdots b_{n+j-1}^2 b_{n+j}^2 \cdots b_{n+j+r-1}^2}
\nonumber\\
&=
\left(\sum_{\bfb'\in \Na^{n-\ell-1}} 
\leb(I_{k-1}(\bfb') ) \right)
\sum_{\substack{b_{k}\cdots b_{k+\ell-1}\geq \vphi(n)\\ b_{n}\cdots b_{n+\ell-1}\geq \vphi(n) }} 
\frac{1}{ b_{k}^2\cdots b_{k+r-1}^2  b_{n}^2 \cdots b_{n+j-1}^2 b_{n+j}^2 \cdots b_{n+j+r-1}^2} \nonumber\\
&=
\sum_{\substack{b_{k}\cdots b_{k+\ell-1}\geq \vphi(n)\\ b_{n}\cdots b_{n+\ell-1}\geq \vphi(n) }} 
\frac{1}{ b_{k}^2\cdots b_{k+r-1}^2  b_{n}^2 \cdots b_{n+j-1}^2 b_{n+j}^2 \cdots b_{n+j+r-1}^2} \\
&\asymp
\frac{\log^{j-1}\vphi(n)}{\vphi(n)} &&\text{(by Theorem \ref{Thm:OverlapEstimatesSeries})}\\
&=
\frac{\log^{k+\ell-n-1}\vphi(n)}{\vphi(n)}.
\end{align*}

\end{enumerate}
 
\end{proof}
\begin{proof}[Proof of Theorem \ref{TEO:MAIN:MEASURE}. Convergence case.]
For any $n\in\Na_{\geq \ell}$ we have
\begin{align*}
\leb(A_n(\vphi) )
&\leq 
\sum_{k=1}^{n-1} \leb(A_{n,k}(\vphi)) \nonumber\\
&=
\sum_{k=1}^{n-\ell} \leb(A_{n,k}(\vphi)) + \sum_{k=n-\ell+1}^{n-1} \leb(A_{n,k}(\vphi)) \nonumber\\
&\asymp
(n-\ell)\frac{\log^{2(\ell-1)}\vphi(n)}{\vphi^2(n)} + \frac{1}{\vphi(n)}\sum_{k=n-\ell+1}^{n-1} \log^{k+\ell-n-1}\vphi(n) \nonumber\\
&\asymp
\frac{n\log^{2(\ell-1)}\vphi(n)}{\vphi^2(n)} + \frac{\log^{\ell-2}\vphi(n)}{\vphi(n)}.
\end{align*}
The desired conclusion follows now from the convergence part of the Borel-Cantelli Lemma.
\end{proof}

\subsection{Divergence case}
 
We now assume the divergence of the series in Theorem \ref{TEO:MAIN:MEASURE}.
We impose two additional hypotheses without losing any generality. First, we suppose that 
\[
\lim_{n\to\infty} \vphi(n)=\infty.
\]
If this were not the case, then every irrational number with unbounded partial quotients would belong to $\clF_{\ell}(\vphi)$ and, thus, $\leb(\clF_{\ell}(\vphi))=1$ by the Borel-Bernstein Theorem (or Birkhoff's Ergodic Theorem). Second, we may assume that \footnote{There seems to be a gap in the proof of this assertion for $\ell=2$ in \cite{TanZhou2024}. Consider $0<\alpha<1$ and $\vphi(n)=n(\log n)^{\alpha}$. Then, $\sum_n \vphi(n)^{-1}=\infty$ and $\vphi(n)< n\log \vphi(n)$ for all large $n$. According to \cite[Section 3.2]{TanZhou2024}, we should consider the function $\psi$ given by $\psi(n)=\max\{\vphi(n), n\log\vphi(n)\}$. However, for all large $n$ we have $n\log \psi(n) = n\log (n\log \vphi(n)) > n\log \vphi(n) = \psi(n)$. This minor gap does not affect the validity of that paper at all.} 
\begin{equation}\label{Eq:WLOG}
\varphi(n) \geq n \log^{\ell-1}\varphi(n)
\quad
\text{ for all }n\in \Na.
\end{equation}
Certainly, for each $n\in\Na$ let $x_n$ be the solution $x$ of the equation
\[
    \frac{x}{\log^{\ell-1}x} = n. 
\]
Then, $(x_n)_{n\geq 1}$ tends to infinity and it is eventually strictly increasing. Let $\psi: \mathbb{N} \to \mathbb{N}$ be given by 
\[
    \psi(n) = \max\{\varphi(n),x_n\}
    \quad
    \text{ for all }n\in\Na.
\]
We can trivially modify at most finitely many terms of $\psi$ to ensure that it is non-decreasing. 
Furthermore, since $\psi \geq \varphi$, we have $\clF_{\ell}(\psi) \subseteq \clF_{\ell}(\varphi)$ and $\psi(n) \to \infty$ as $n\to\infty$. It is also clear that $\psi(n) \geq n\log^{\ell-1} \psi(n)$ for all $n\in\Na$, because $x \mapsto \frac{x}{\log^{\ell-1}x}$ is strictly increasing for large $x$. Lastly, let us show that
\[
    \sum_{n=1}^\infty \left(n \left(\frac{\log^{\ell-1}\psi(n)}{\psi(n)}\right)^2 + \frac{\log^{\ell-2}\psi(n)}{\psi(n)} \right)= \infty.
\]
The divergence is obvious if $\varphi(n)\geq x_n$ for all but finitely many $n \in \mathbb{N}$. If this is not the case, then there is a strictly increasing sequence of natural numbers $(m_k)_{k\geq 1}$ such that
\[
m_{k+1}>2m_k \text{ and } \varphi(m_k) < x_{m_k}
\;\text{ for all }\;k\in\Na.
\]
Hence, for large $K=K(\ell)\in\Na$, we have
\begin{align*}
    \sum_{n=1}^\infty \left(n \left(\frac{\log^{\ell-1}\psi(n)}{\psi(n)}\right)^2 + \frac{\log^{\ell-2}\psi(n)}{\psi(n)} \right)
    &\geq 
    \sum_{k=K}^\infty 
    \sum_{n=\left[\frac{m_{k+1}}{2}\right]+1}^{ m_{k+1}} 
     n \left(\frac{\log^{\ell-1}\psi(n)}{\psi(n)}\right)^2\\
    &\geq 
    \sum_{k=K}^\infty \frac{m_{k+1}}{2} 
    \sum_{\left[\frac{m_{k+1}}{2}\right]+1 \leq n \leq m_{k+1}} 
    \left(\frac{\log^{\ell-1}\psi(n)}{\psi(n)}\right)^2\\
    &
    \geq 
    \sum_{k=K}^\infty \frac{m_{k+1}}{2} 
    \sum_{n=\left[\frac{m_{k+1}}{2}\right]+1}^{ m_{k+1}} 
    \left(\frac{\log^{\ell-1}x_{m_{k+1}}}{x_{m_{k+1}}}\right)^2\\
    &= 
    \sum_{k=K}^\infty \frac{m_{k+1}}{2} 
    \sum_{n=\left[\frac{m_{k+1}}{2}\right]+1}^{ m_{k+1}} 
    \left(\frac{1}{m_{k+1}}\right)^2 \\
    &\geq 
    \sum_{k=K}^\infty \frac{m_{k+1}}{2} 
        \left(\frac{1}{m_{k+1}}\right)^2
        \left( \frac{m_{k+1}}{2} - 1\right)\\
    &=
    \sum_{k=K}^\infty 
    \frac{1}{2m_{k+1}}
        \left( \frac{m_{k+1}}{2} - 1\right)\\
    &=
    \sum_{k=K}^\infty 
    \frac{1}{4}
        \left(1 - \frac{2}{m_{k+1}}  \right)= \infty.
\end{align*}
\subsubsection{Lemmas for the divergence case}
Take $t\in\Na$, $\bfb=(b_1,\ldots, b_t)\in \Na^t$, $I_t=I_t(\bfb)$, and let $M=M(\bfb)\in\Na_{\geq t+1}$ be such that $\vphi(M)>\left(\max\{b_1,\ldots, b_t\}\right)^{\ell}$.
\begin{lem01}\label{Lem:Divergence:01}
If $n\in\Na_{\geq M}$, then
\[
\leb\left( A_n\cap I_t\right)
\gg_{\ell}
\leb(I_t)
\left( \frac{n\log^{2(\ell-1)}\vphi(n)}{\vphi^2(n)} +  \frac{\log^{\ell-2}\vphi(n)}{\vphi(n)} \right).
\]
\end{lem01}
\begin{proof}
Take $n\in\Na_{\geq M}$. Define the set
\[
B_{n}
:=
\left\{x\in I_t: 
\begin{matrix}
a_{n-1}(x)\cdots a_{n+\ell-2}(x) \geq \vphi(n), \\ 
a_{n}(x)\cdots a_{n+\ell-1}(x) \geq \vphi(n)
\end{matrix}
\right\}.
\]
For each $i\in\Na$ such that $M\leq i<n-1$ and $i\equiv n-1 \pmod \ell$, let $B_i$ be the set of those $x\in I_t$ satisfying the next properties:
\begin{enumerate}[1.]
\item $a_{i}(x)\cdots a_{i+\ell-1}(x)\geq \vphi(n)$, 
\item For each $k\in\Na$ such that $0< k < \frac{n-1-i}{\ell}$ we have $a_{i+k\ell }(x)\cdots a_{i + (k+1)\ell-1}(x)< \vphi(n)$.
 \item $a_{n-1}(x)\cdots a_{n+\ell-2}(x)<\vphi(n)$ and $a_{n}(x)\cdots a_{n+\ell-1}(x)\geq \vphi(n)$.
\end{enumerate}
That is, we impose restrictions on non-overlapping contiguous blocks of length $\ell$. Clearly, we have 
\begin{equation}\label{Eq:Conten}
B_n\cup \bigcup_{\substack{M\leq  i<n-1\\ i\equiv n-1 (\ell)}} B_i
\subseteq
A_n\cap I_t,
\end{equation}
and the union is pairwise disjoint. 
First, Theorem \ref{Thm:OverlapEstimatesSeries} yields
\begin{equation}\label{Eq:BndBn}
\leb(B_{n})
\asymp_{\ell}
\frac{\log^{\ell-2} \vphi(n) }{\vphi(n)}
\leb(I_t).
\end{equation}
For the measure of the remaining terms, put 
\[
P = \left[ \frac{n-1-M}{\ell}\right] 
\quad\text{ and }\quad 
r_1 = n-2-P\ell,
\]
so we are dealing with $P+1$ disjoint contiguous blocks of length $\ell$ with the last one of them being $a_{n-1},\ldots, a_{n+\ell-1}$.
Given $i\in \Na$ with $M\leq i<n-1$ and $i\equiv n-1\pmod \ell$, we define $k_i=\frac{n-1-i}{\ell}$.
Then, letting $\bfc_0, \ldots, \bfc_{k_i-1}$ represent elements in $\Na^{\ell}$ in the forthcoming expressions and writing 
\[
\Pi(\bfx) = \prod_{j=1}^{\ell} x_j 
\quad\text{ for all }\quad 
\bfx=(x_1, \ldots, x_\ell)\in \Na^{\ell},
\]
we have
\[
B_i
=
\bigcup_{\bfe\in\Na^{r_1 + (P-k)\ell }}
\bigcup_{\prod(\bfc_{0})\geq \vphi(n)}
\bigcup_{\substack{\Pi(\bfc_{j})< \vphi(n) \\ 1\leq j\leq k_i-1  }}
\bigcup_{\substack{a_{n-1}\cdots a_{n+\ell-2}<\vphi(n) \\ a_{n}\cdots a_{n+\ell-1}\geq \vphi(n) } }
I_{n+\ell-1}(\bfb,\bfe, \bfc_0, \cdots, \bfc_{k_i-1},a_{n-1},a_{n},\cdots ,a_{n-\ell-2,}a_{n-\ell-1}).
\]
There are no actual restrictions on $\bfe$ other than its length, then $\leb (B_i)$ is, up to multiplicative constants, equivalent to
\[
\leb(I_{t}(\bfb))
\sum_{\prod(\bfc_{0})\geq \vphi(n)}
\sum_{\substack{\Pi(\bfc_{j})< \vphi(n) \\ 1\leq j\leq k_i-1  }}
\sum_{\substack{a_{n-1}\cdots a_{n+\ell-2}<\vphi(n) \\ a_{n}\cdots a_{n+\ell-1}\geq \vphi(n) } }
\leb(I_{(k_i+1)\ell+1}(\bfc_0, \cdots, \bfc_{k_i-1},a_{n-1},a_{n},\cdots ,a_{n-\ell-2,}a_{n-\ell-1})).
\]
By Lemma \ref{Lem:HWX}, Lemma \ref{Lem:SerEstNu}, and Lemma \ref{Lem:CF02}, we have
\[
\leb(B_i)
\asymp
\leb(I_{t}(\bfb))
\frac{\log^{2(\ell-1)}\vphi(n)}{\vphi^2(n)}
\sum_{\substack{\Pi(\bfc_{j})< \vphi(n) \\ 1\leq j\leq k_i-1  }}
\leb(I_{n+1}\left( \bfc_1\cdots \bfc_{k_i-1 })\right).
\]
To estimate the last sum, note that there is a positive constant $\kappa=\kappa(\ell)>0$ such that, for all $m\in\Na$ and $\bff\in \Na^m$,
\begin{align*}
\sum_{\substack{\bfa\in\Na^{\ell} \\ \prod(\bfa)< \vphi(n) }}
\leb(I_{m+\ell} (\bff \bfa))
&= 
\leb(I_m(\bff))
- 
\sum_{\substack{\bfa\in\Na^{\ell} \\ \prod(\bfa)\geq \vphi(n) }}
\leb(I_{m+\ell} (\bff \bfa)) \\
&\geq 
\leb(I_m(\bff))
\left(1 - 
\kappa\frac{\log^{\ell-1}\vphi(n)}{\vphi(n)}\right) \\
&\geq 
\leb(I_m(\bff))
\left( 1 -  \frac{\kappa}{n}\right).
\end{align*}
Therefore, assuming $n\in\Na$ is large enough,
\begin{align*}
\leb(B_i)
&\gg
\leb(I_{t}(\bfb))
\frac{\log^{2(\ell-1)}\vphi(n)}{\vphi^2(n)}
\left( 1 - \frac{\kappa}{n}\right)^{k_i-1} \\
&=
\leb(I_{t}(\bfb))
\frac{\log^{2(\ell-1)}\vphi(n)}{\vphi^2(n)}
\left( 1 - \frac{\kappa}{n}\right)^{\frac{n-1-i}{\ell} - 1} \\
&\geq
\leb(I_{t}(\bfb))
\frac{\log^{2(\ell-1)}\vphi(n)}{\vphi^2(n)}
\left( 1 - \frac{\kappa}{n}\right)^{ n-i}.
\end{align*}

Using elementary identities about geometric sums and $0<e^{-\kappa}<1$, we conclude 
\begin{equation}\label{Eq:BndBn02}
\leb\left(\bigcup_{\substack{M\leq i<n-1\\ i\equiv n-1 (\ell)}} B_i\right)
\gg
\leb(I_t)
\frac{n\log^{2(\ell-1)}\vphi(n)}{\vphi^2(n)}.
\end{equation}
The lemma now follows from \eqref{Eq:Conten}, \eqref{Eq:BndBn}, and \eqref{Eq:BndBn02}.
\end{proof} 

\begin{lem01}\label{Lem:Divergence:02}
For large $N\in\Na$ we have
\begin{equation}\label{Eq:Lem:4.3}
\sum_{M\leq m< n \leq N} \leb\left(A_m\cap A_n\cap I_t\right)
\ll_{\ell}
\sum_{n=M}^N \leb\left( A_n\cap I_t \right)
+
\leb(I_t)^{-1} \left(\sum_{n=M}^N \leb\left( A_n\cap I_t \right) \right)^2.
\end{equation}
\end{lem01}
\begin{proof}
We partition the index set 
\[
\clI:=\{(m,n)\in\Na^2: M\leq m< n \leq N\}
\]
into
\[
\clI_1:=\{(m,n)\in \clI: M\leq m \leq n - (2\ell-1)\} 
\;\text{ and }\;
\clI_2:=\{(m,n)\in \clI: n-2\ell+2\leq m \leq n-1 \}  \},
\]
and write
\[
S_1:=\sum_{(m,n)\in \clI_1} \leb\left(A_m\cap A_n\cap I_t\right)
\;\text{ and }\;
S_2:=\sum_{(m,n)\in \clI_2} \leb\left(A_m\cap A_n\cap I_t\right).
\]
It is clear that 
\begin{equation}\label{Eq:Div:Lem:02:01}
S_2 \leq  2\ell\sum_{n=M}^N \leb\left(A_n\cap I_t\right).    
\end{equation}

Now we give an upper bound for $S_1$. 
For each $(m,n)\in \clI_1$, define the sets $U_1(m,n)$, $U_2(m,n)$, $U_3(m,n)$ as follows:
\begin{enumerate}[1.]
\item Let $U_1(m,n)$ be the set of irrational numbers $x\in I_t$ such that
\begin{itemize}
\item[•] $a_m(x)\cdots a_{m+\ell-1}(x)\geq \vphi(m)$,
\item[•] For some $k<m$ we have $a_{k}(x)\cdots a_{k+\ell-1}(x)\geq \vphi(n)$,
\item[•] $a_n(x)\cdots a_{n+\ell-1}(x)\geq \vphi(n)$.
\end{itemize}

\item Let $U_2(m,n)$ be the set of irrational numbers $x\in I_t$ such that
\begin{itemize}
\item[•] For some $j<m$ we have $a_{j}(x)\cdots a_{j+\ell-1}(x)\geq \vphi(m)$,
\item[•] $a_m(x)\cdots a_{m+\ell-1}(x)\geq \vphi(n)$,
\item[•] $a_n(x)\cdots a_{n+\ell-1}(x)\geq \vphi(n)$.
\end{itemize}

\item Let $U_3(m,n)$ be the set of irrational numbers $x \in I_t$ such that
\begin{itemize}
\item[•] For some $j<m$ we have $a_{j}(x)\cdots a_{j+\ell-1}(x)\geq \vphi(m)$,
\item[•] $a_m(x)\cdots a_{m+\ell-1}(x)\geq \vphi(m)$,
\item[•] For some $m<k<n$ we have $a_{k}(x)\cdots a_{k+\ell-1}(x)\geq \vphi(n)$,
\item[•] $a_n(x)\cdots a_{n+\ell-1}(x)\geq \vphi(n)$.
\end{itemize}
\end{enumerate}
As a consequence, 
\[
A_m\cap A_n\cap I_t
\subseteq 
U_1(m,n)\cup U_2(m,n) \cup U_3(m,n).
\]
We claim that
\begin{equation}\label{Eq:EstU1}
\leb(U_1(m,n))
\ll
\leb(I_t) \frac{\log^{2(\ell-1)}\vphi(n)}{\vphi^2(n)},    
\end{equation}
\begin{equation}\label{Eq:EstU2}
\leb(U_2(m,n))
\ll
\leb(I_t) \frac{\log^{2(\ell-1)}\vphi(n)}{\vphi^2(n)}, \;\text{ and }
\end{equation}
\begin{align}\label{Eq:EstU3}
\leb(U_3(m,n))
&\ll
\leb(I_t)
\frac{\log^{2(\ell-1)}\vphi(n)}{\vphi^2(n)}  +\nonumber\\
&\quad\; + 
\leb(I_t)\left( \frac{n\log^{2(\ell-1)}  \vphi(n)}{\vphi^2(n)} + \frac{\log^{\ell-2}\vphi(n)}{\vphi(n)} \right)
\left( \frac{m\log^{2(\ell-1)}  \vphi(m)}{\vphi^2(m)} + \frac{\log^{\ell-2}\vphi(m)}{\vphi(m)} \right).
\end{align}

First, we show \eqref{Eq:EstU1}. 
When $M\leq k \leq m-\ell$, we apply \eqref{Eq:Series:01} to get
\[
\sum_{\substack{a_k \cdots a_{k+ \ell- 1}\geq \vphi(n) \\ a_m \cdots a_{m+ \ell- 1}\geq \vphi(m) \\ a_n\cdots a_{n + \ell- 1}\geq \vphi(n)} }
\frac{1}{a_k^2\cdots a_{k+\ell-1}^2 a_m^2\cdots a_{m+\ell-1}^2 a_{n}^2\cdots a_{n+\ell-1}^2} 
\asymp
\frac{\log^{\ell-1}\vphi(m)}{\vphi(m)}
\frac{\log^{2(\ell-1)}\vphi(n)}{\vphi^2(n)}.
\]
When $m-\ell+1\leq k\leq m-1$, we have
\begin{align*}
\sum_{\substack{ a_{k}\cdots a_{k+\ell-1}\geq \vphi(n)\\ a_m\cdots a_{m+\ell-1}\geq \vphi(m)\\a_n\cdots a_{n+\ell-1}\geq \vphi(n)}}
\frac{1}{(a_{k}^2 \cdots a_{m + \ell-1}^2) a_{n}^2\cdots a_{n+\ell-1}^2} 
&\leq
\sum_{\substack{ a_{k}\cdots a_{k+\ell-1}\geq \vphi(n)\\ a_n\cdots a_{n+\ell-1}\geq \vphi(n)}}
\frac{1}{(a_{k}^2 \cdots a_{m + \ell-1}^2) a_{n}^2\cdots a_{n+\ell-1}^2} \\
&\ll
\frac{\log^{2(\ell-1)}\vphi(n)}{\vphi^2(n)}
\sum_{a_{k+\ell},\ldots, a_{m+\ell-1}\in\Na}
\frac{1}{a_{k+\ell}^2 \cdots a_{m + \ell-1}^2} \\
&\asymp
\frac{\log^{2(\ell-1)}\vphi(n)}{\vphi^2(n)}.
\end{align*}
The above inequalities along with \eqref{Eq:WLOG} yield
\begin{align*}
\leb(U_1(m,n))
&\ll 
m\frac{\log^{\ell-1}\vphi(m)}{\vphi(m)}
\frac{\log^{2(\ell-1)}\vphi(n)}{\vphi^2(n)}
+
\ell \frac{\log^{2(\ell-1)}\vphi(n)}{\vphi^2(n)} \\
&\ll 
\frac{\log^{2(\ell-1)}\vphi(n)}{\vphi^2(n)},
\end{align*}
which is \eqref{Eq:EstU1}. 
We have used that for the considered values of $k$ considered, the blocks $(a_k,\ldots, a_{k+\ell-1})$ and $(a_n,\ldots, a_{n+\ell-1})$ do not overlap. 

We work with two different series to estimate $\leb(U_2(m,n))$. 
First, when $M\leq j\leq m- \ell$, the blocks $(a_j,\ldots, a_{j+\ell-1})$ and $(a_m,\ldots, a_{m+\ell-1})$ do not overlap and
\begin{align*}
\sum_{j=M}^{m-\ell} 
\sum_{\substack{ a_j\cdots a_{j+\ell-1}\geq \vphi(m) \\ a_m\cdots a_{m+\ell-1}\geq \vphi(n)}}
\frac{1}{a_j^2 \cdots a_{j+\ell-1}^2 a_m^2\cdots a_{m+\ell-1}^2}
&\asymp
m \frac{\log^{\ell-1}\vphi(m)}{\vphi(m)}\frac{\log^{\ell-1}\vphi(n)}{\vphi(n)} \\
&\leq
\frac{\log^{\ell-1}\vphi(n)}{\vphi(n)} &&\text{(by \eqref{Eq:WLOG}}).
\end{align*}
In the forthcoming series, the product $a_j^2 \cdots a_{m+\ell-1}^2$: if $j+\ell-1\leq m$, then 
\[
a_j^2 \cdots a_{m+\ell-1}^2
=
(a_j\cdots a_{j+\ell-1})(a_m \cdots a_{m+\ell-1});
\]
if $j+\ell-1> m$, then $a_j^2 \cdots a_{m+\ell-1}^2$ is the product of the block $(a_{j},a_{j+1}, \ldots, a_{m+\ell-1})$. 
We use the same convention for other combinations of the parameters $j, k, m, n$.

When $m- \ell +1 \leq j \leq m-1$, the blocks $(a_j,\ldots, a_{j+\ell-1})$ and $(a_m,\ldots, a_{m+\ell-1})$ overlap and
\begin{align*}
\sum_{j=m - \ell +1}^{m-1} 
\sum_{\substack{ a_j\cdots a_{j+\ell-1}\geq \vphi(m) \\ a_m\cdots a_{m+\ell-1}\geq \vphi(n)}}
\frac{1}{a_j^2 \cdots a_{m+\ell-1}^2}
&\leq 
\sum_{j=m - \ell +1}^{m-1} 
\sum_{\substack{ a_j,\ldots, a_{j+\ell-1}\in\Na  \\ a_m\cdots a_{m+\ell-1}\geq \vphi(n)}}
\frac{1}{a_j^2 \cdots a_{m+\ell-1}^2}\\
&\ll_{\ell}
\frac{\log^{\ell-1}\vphi(n)}{\vphi(n)}.
\end{align*}
Finally, applying Lemma \ref{Lem:HWX}, we conclude \eqref{Eq:EstU2}: 
\begin{align*} 
\leb &(U_2(m,n)) \\
&\ll 
\left(
\sum_{ a_n\cdots a_{n+\ell-1}\geq \vphi(n)}
\frac{1}{a_n^2 \cdots a_{n+\ell-1}^2}
\right)
\left(
\sum_{j=M}^{m-\ell} 
\sum_{\substack{ a_j\cdots a_{j+\ell-1}\geq \vphi(m) \\ a_m\cdots a_{m+\ell-1}\geq \vphi(n)}}
\frac{1}{a_j^2 \cdots a_{m+\ell-1}^2}
+
\sum_{j=m - \ell +1}^{m-1} 
\sum_{\substack{ a_j\cdots a_{j+\ell-1}\geq \vphi(m) \\ a_m\cdots a_{m+\ell-1}\geq \vphi(n)}}
\frac{1}{a_j^2 \cdots a_{m+\ell-1}^2}
\right) \\
&\ll 
\frac{\log^{2(\ell-1)}\vphi(n)}{\vphi^2(n)}.
\end{align*}

Next, we show \eqref{Eq:EstU3}. First, we consider the values of $k$ for which the block $(a_k,\ldots, a_{k+\ell-1})$ does not overlap neither $(a_n,\ldots, a_{n+\ell-1})$ nor $(a_m,\ldots, a_{m+\ell-1})$. In other words, we assume $m+\ell \leq k \leq n-\ell$. This leads us to the series
\[
\sum_{j=M}^{m-1}
\sum_{k=m+\ell}^{n-\ell}
\sum_{\substack{a_j\cdots a_{j+\ell-1}\geq \vphi(m)\\a_m\cdots a_{m+\ell-1}\geq \vphi(m)\\a_k\cdots a_{k+\ell-1}\geq \vphi(n)\\a_n\cdots a_{n+\ell-1}\geq \vphi(n)}}
\frac{1}{a_j^2 \cdots a_{m+\ell-1}^2a_k^2\cdots a_{k+\ell-1}^2a_n^2\cdots a_{n+\ell-1}^2},
\]
which is asymptotically equivalent to
\begin{align}
(n-m)
\left( \frac{\log^{\ell-1}\vphi(n)}{\vphi(n)}\right)^2
&\sum_{j=M}^{m-1}
\sum_{\substack{a_j\cdots a_{j+\ell-1}\geq \vphi(m)\\a_m\cdots a_{m+\ell-1}\geq \vphi(m)}}
\frac{1}{a_j^2 \cdots a_{m+\ell -1}^2} 
\nonumber\\
&\ll
n
\frac{\log^{2(\ell-1)}\vphi(n)}{\vphi^2(n)}
\left( m \frac{\log^{2(\ell-1)}\vphi(m)}{\vphi(m)} + \frac{\log^{\ell-2}\vphi(m)}{\vphi(m)}\right). \label{Eq:U3:a}
\end{align}
The last estimate is obtained as in the convergence case. 

Now suppose that $n-\ell+1\leq k \leq n-1$. This means that $(a_k,\ldots, a_{k+\ell-1})$ and $(a_n,\ldots, a_{n+\ell-1})$ overlap but $(a_k,\ldots, a_{k+\ell-1})$ and $(a_m,\ldots, a_{m+\ell-1})$ do not (by the definition of $\clI_1$), hence
\begin{align}
\sum_{j=M}^{m-1}
&\sum_{k=n-\ell+1}^{n-1}
\sum_{\substack{a_j\cdots a_{j+\ell-1}\geq \vphi(m)\\a_m\cdots a_{m+\ell-1}\geq \vphi(m)\\a_k\cdots a_{k+\ell-1}\geq \vphi(n)\\a_n\cdots a_{n+\ell-1}\geq \vphi(n)}}
\frac{1}{a_j^2 \cdots a_{m+\ell-1}^2a_k^2\cdots a_{n+\ell-1}^2} \ll \nonumber\\
&\ll
\frac{\log^{\ell-2}\vphi(n)}{\vphi(n)}
\sum_{j=M}^{m-1}
\sum_{\substack{a_j\cdots a_{j+\ell-1}\geq \vphi(m) \\ a_m\cdots a_{m+\ell-1}\geq \vphi(m) }}
\frac{1}{a_j^2 \cdots a_{m+\ell-1}^2} \nonumber\\
&\ll
\frac{\log^{\ell-2}\vphi(n)}{\vphi(n)}
\left( m \frac{\log^{2(\ell-1)}\vphi(m)}{\vphi(m)} + \frac{\log^{\ell-2}\vphi(m)}{\vphi(m)}\right). \label{Eq:U3:b}
\end{align}
For the first estimate, we used that the length overlap of $(a_k, \ldots, a_{k+\ell-1})$ and $(a_n, \ldots, a_{n+\ell-1})$ is at most $\ell-1$.
The last estimate is obtained as in the case $m+\ell \leq k \leq n-\ell$.

Lastly, we deal with the values of $k$ for which the blocks $(a_m,\ldots, a_{m+\ell-1})$ and $(a_k, \ldots, a_{k+\ell-1})$ overlap. This happens if and only if $k\leq m+\ell-1$. In this case, 
\[
k+ \ell-1\leq m+2(\ell-1)<n,
\]
hence $(a_k, \ldots, a_{k+\ell-1})$ and $(a_n, \ldots, a_{n+\ell-1})$ do not overlap. For any $j\in \{M, \ldots, m-1\}$, $k\in \{m+1, \ldots, m+\ell-1\}$ define
\[
D_{j,k}
:=
\left\{
x\in I_t:
\begin{matrix}
a_j(x)\cdots a_{j+\ell-1}(x)\geq \vphi(m), & a_m(x)\cdots a_{m+\ell-1}(x)\geq \vphi(m), \\
a_k(x)\cdots a_{k+\ell-1}(x)\geq \vphi(n), & a_n(x)\cdots a_{n+\ell-1}(x)\geq \vphi(n)
\end{matrix}
\right\}.
\]
Then, we have
\[
\bigcup_{j=M}^{m-1} D_{j,k}
\subseteq
\left\{
x\in I_t:
\begin{matrix}
a_k(x)\cdots a_{k+\ell-1}(x)\geq \vphi(n),\\  
a_n(x)\cdots a_{n+\ell-1}(x)\geq \vphi(n)
\end{matrix}
\right\}
=A_{n,k}\cap I_t.
\]
We conclude
\begin{equation}\label{Eq:U3:c}
\leb
\left( \bigcup_{k=m+1}^{m+\ell-1}\bigcup_{j=M}^{m-1} D_{j,k}\right)
\leq 
\sum_{k=m+1}^{m+\ell-1} \leb(A_{n,k}\cap I_t)
\asymp_{\ell}
\frac{\log^{2(\ell-1)}\vphi(n)}{\vphi^2(n)}\leb(I_t).    
\end{equation}
Then, \eqref{Eq:EstU3} follows from \eqref{Eq:U3:a}, \eqref{Eq:U3:b}, and \eqref{Eq:U3:c}.

Putting \eqref{Eq:EstU1}, \eqref{Eq:EstU2}, and \eqref{Eq:EstU3} together, we have
\begin{align*}
S_1
&\leq 
\sum_{(m,n)\in \clI_1} \leb(U_1(m,n)) +\leb(U_2(m,n))+ \leb(U_3(m,n)) \\
&\ll
\leb(I_t)
\sum_{(m,n)\in \clI_1}
\frac{\log^{2(\ell-1)}\vphi(n)}{\vphi^2(n)}  +  \left( \frac{n\log^{2(\ell-1)}  \vphi(n)}{\vphi^2(n)} + \frac{\log^{\ell-2}\vphi(n)}{\vphi(n)} \right)
\left( \frac{m\log^{2(\ell-1)}  \vphi(m)}{\vphi^2(m)} + \frac{\log^{\ell-2}\vphi(n)}{\vphi(m)}\right).
\end{align*}
Also, by the definition of $\clI_1$ and Lemma \ref{Lem:Divergence:01}, we have
\begin{align*}
\sum_{(m,n)\in \clI_1}&
\left(\frac{n\log^{2(\ell-1)}  \vphi(n)}{\vphi^2(n)} + \frac{\log^{\ell-2}\vphi(n)}{\vphi(n)} \right)
\left( \frac{m\log^{2(\ell-1)}  \vphi(m)}{\vphi^2(m)} + \frac{\log^{\ell-2}\vphi(n)}{\vphi(m)}\right) \\
&=
\sum_{n=M + 2\ell-1}^{N}
\sum_{m=M}^{n-(2\ell-1)}
\left(\frac{n\log^{2(\ell-1)}  \vphi(n)}{\vphi^2(n)} + \frac{\log^{\ell-2}\vphi(n)}{\vphi(n)} \right)
\left( \frac{m\log^{2(\ell-1)}  \vphi(m)}{\vphi^2(m)} + \frac{\log^{\ell-2}\vphi(n)}{\vphi(m)}\right)\\
&\ll
\frac{1}{\leb(I_t)^2}
\left( \sum_{n=M}^N \leb(A_n\cap I_t)\right)^2.
\end{align*}
We also deduce from Lemma \ref{Lem:Divergence:01} that 
\[
\sum_{(m,n)\in \clI_1}
\frac{\log^{2(\ell-1)}\vphi(n)}{\vphi^2(n)}
\ll
\sum_{n=M}^N
\frac{n\log^{2(\ell-1)}\vphi(n)}{\vphi^2(n)}
\ll
\frac{1}{\leb(I_t)}
\sum_{n=M}^N \leb(A_n\cap I_t).
\]
Therefore, we have
\begin{equation}\label{Eq:Div:Lem:02:02}
S_1\ll \sum_{n=M}^N \leb(A_n\cap I_t) + \frac{1}{\leb(I_t)}\left( \sum_{n=M}^N \leb(A_n\cap I_t)\right)^2.
\end{equation}
We combine \eqref{Eq:Div:Lem:02:01} and \eqref{Eq:Div:Lem:02:02}, and the result follows.
\end{proof}

\begin{proof}[Proof of Theorem \ref{TEO:MAIN:MEASURE}. Divergence case]
Take any $t\in\Na$, any $\bfb\in\Na^t$, and consider $I_t=I_t(\bfb)$ and $M$ as above. 
For each $n\in\Na$, put $E_n:=A_{M-1+n} \cap I_t$. 
By Lemma \ref{Lem:Divergence:01} and the divergence assumption, we have $\displaystyle\sum_{n\geq 1} \leb(E_n)=\infty$;
therefore, Lemma \ref{Lem:Divergence:02} and the Chung-Erd\H{o}s inequality (Lemma \ref{Le:ErdosChung}) yield
\[
\leb\left( \limsup_{n\to\infty} E_n\right)
\gg_{\ell}
\leb(I_t).
\]
Then, we apply Lemma \ref{Le:Knopp} to the collection $\clC$ is the collection of fundamental intervals to conclude $\leb(\clF_{\ell}(\vphi))=1$.

\end{proof}

\section{Proof of Theorem \ref{TEO:MAIN:HD}}\label{Section:ProofHD}
Theorem \ref{TEO:MAIN:HD} is obvious when $B_{\vphi}=1$ or $B_{\vphi}=\infty$ in view of $\clF_2(\vphi)\subseteq \clF_3(\vphi)\subseteq \clE_3(\vphi)$ and theorems \ref{TEO:HWX:HD} and \ref{TEO:TZ:HD}.
We, thus, only deal with the case $1<B_{\vphi}<\infty$.
We further restrict our calculations to functions of the form $\vphi_B:\Na\to\RE$, $\vphi_B(n)=B^n$ for a fixed real number $B>1$. 
Standard arguments (see cf. \cite{HuangWuXu2020, WangWu2008}) allow us to conclude the general case.

\begin{lem01}\label{Lem:HD:Main}
Let $B>1$ be a real number, then 
\[
\hdim \clF_3(\vphi_B)
= 
\inf\left\{s\geq 0: P(T,-s\log|T'| - g(s)\log B)\leq 0\right\}.
\]
\end{lem01}

\subsection{Lower bound}
To compute the lower bound, we construct a family of subsets $S(x,y)$ of $\clF_3(\vphi_B)$ parametrized by certain pairs $(x,y)\in\RE^2$. 
We then maximize $(x,y)\mapsto \hdim S(x,y)$. 
For any $x,y\in \RE$ such that $0\leq x\leq 1$, $0\leq y\leq 1$, and $0\leq x+y\leq 1$, define 
\[
A_0(x,y) = B^{1-x-y}, \quad
A_1(x,y) = B^{x}, \quad
A_2(x,y) = B^{y}, \quad
A_3(x,y) = B^{1-x-y}.
\]
The set 
\[
S(x,y)
=
S_4(A_0(x,y),A_1(x,y),A_2(x,y),A_3(x,y))
\]
defined in \eqref{Eq:Def:SetS} is, thus, contained in $\clF_3(\vphi_B)$; so $\hdim S(x,y) \leq \hdim \clF_3(B)$; therefore,
\[
\sup_{ \substack{0\leq x,y\leq 1 \\ 0\leq x+y\leq 1 } }\hdim S(x,y)\leq  \hdim \clF_3(\vphi_B).
\]
To compute this supremum, for each $\frac{1}{2}<s<1$ and $(x,y)$ such that $x,y$ are non-negative and $x+y\leq 1$, define the potentials
\begin{align*}
f_{0,s}(x,y) 
&= 
-s(1-x-y) \log B -s \log|T'|, \\
f_{1,s}(x,y) 
&= 
-s(1-y) \log B + (1-s)(1-x-y) \log B - s\log|T'|, \\
f_{2,s}(x,y) 
&= 
-s \log B + (1-s)(1-y) \log B - s\log|T'|, \\
f_{3,s}(x,y) 
&= -s(2-x-y) \log B +(1-s)\log B -s \log|T'|,
\end{align*}
and, for $i\in \{0,1,2,3\}$,
\[
d_i(x,y)
=
\inf\left\{s\geq 0: P(T, f_{i,s}(x,y))\leq 0\right\}.
\]
Then, by Lemma \ref{Le:MumtazNikita}, we must compute
\[
\sup_{ \substack{0\leq x,y\leq 1\\ 0\leq x+y\leq 1 } }
\min_{0\leq i\leq 3} d_i(x,y).
\]
To simplify the problem, rather than $f_{j,s}$, we consider the functions
\begin{align*}
F_{0,s}(x,y) 
&= 
-s(1-x-y), \\
F_{1,s}(x,y) 
&= 
-s(1-y) + (1-s)(1-x-y) 
= -2s+1+ (s-1)x + (2s-1)y,  \\
F_{2,s}(x,y) 
&= 
-s   + (1-s)(1-y) 
=
- 2 s + 1 +(s-1)y, \\
F_{3,s}(x,y) 
&= -s(2-x-y) +(1-s) 
=
- 2 s + 1 -s(1-x-y).
\end{align*}
Therefore, $f_{j,s}(x,y) = F_{j,s}(x,y)\log B - s\log |T'|$ for $j\in \{0,1,2,3\}$.
Since the potentials are increasing in $F_{j,s}(x,y)$, our problem is to find
\[
\sup_{ \substack{0\leq x,y\leq 1 \\ 0\leq x+y\leq 1 } }
\min\{ F_{0,s}(x,y), F_{1,s}(x,y), F_{2,s}(x,y), F_{3,s}(x,y)\}.
\]
Clearly, $F_{0,s}(x,y)\geq  F_{3,s}(x,y)$ for all admissible values of $(x,y)$ and $\frac{1}{2}\leq s\leq 1$, so our goal is to calculate
\[
\sup_{ \substack{0\leq x,y\leq 1\\ 0\leq x+y\leq 1 } }
\min\{F_{1,s}(x,y), F_{2,s}(x,y), F_{3,s}(x,y)\}.
\]
Long yet completely elementary computations allow us to determine the regions where each of the three functions is the minimum. 
In fact, we divide the triangle in the $xy$-plane determined by $0\leq x,y\leq 1$ and $0\leq x+y\leq 1$ into three or two regions depending on $s$.
From this and the signs of the coefficients of $x$ and $y$ on each $F_{j,s}$, we conclude that the supremum is in fact a maximum and that it is attained where all the functions in the minimum coincide. 
This maximum is 
\[
-\frac{3s^3-5s^2+4s-1}{s^2-s+1}=-g(s)
\]
and it is achieved when
\begin{equation}\label{Eq:DefXY}
x=\frac{s^2}{s^2-s+1}
\quad\text{ and }\quad 
y=\frac{s-s^2}{s^2-s+1}.    
\end{equation}
Therefore,
\[
\inf\{s\geq 0: P(T, -g(s)\log B -s \log |T'|)\leq 0\}
\leq 
\hdim \clF_3(\vphi_B).
\]

\subsection{Upper bound.}
For clarity, we start by recalling part of the proof of the upper bound in \cite[Theorem 1.7]{HuangWuXu2020}.
Take a real number $A>1$.
For $n\in\Na$, define
\[
E_2(\vphi_A;n) = \left\{x\in [0,1):A^n\leq a_n(x)a_{n+1}(x)\right\},
\]
and, for $1<\beta\leq A$,
\begin{align*}
E_1(\vphi_{\beta};n) &= \left\{x\in [0,1):\beta^n\leq a_n(x)\right\}, \\
F_2(\vphi_{\beta}, \vphi_A;n) &= \left\{x\in [0,1):1\leq a_n(x)<\beta^n, \, A^n \leq a_{n}(x)a_{n+1}(x)\right\}.
\end{align*}
Note that $\clE_1(\vphi_{\beta}) = \displaystyle\limsup_{n\to\infty} E_1(\vphi_{\beta};n)$; then, writing
\[
F_2(\vphi_{\beta}, \vphi_A) = \limsup_{n\to\infty} F_2(\vphi_{\beta}, \vphi_A;n),
\]
we have 
\[
\clE_2(\vphi_A)\subseteq E_1(\vphi_{\beta})\cup F_2(\vphi_{\beta}, \vphi_A).
\]

For $n\in\Na_{\geq 2}$, and $\bfa\in\Na^{n-1}$, the set $E_1(\vphi_{\beta};n)\cap I_{n-1}(\bfa)$ is a single interval and
\[
\left|E_1(\vphi_{\beta};n)\cap I_{n-1}(\bfa)\right|
=
\left|\bigcup_{a_n\geq \beta^n} I_{n}(\bfa a_n)\right|
\asymp
\frac{1}{q_{n-1}^{2}(\bfa) \beta^n}.
\]
The set  $F_2(\vphi_{\beta}, \vphi_A; n)\cap I_{n-1}(\bfa)$ is a union of separated intervals; more precisely,  
\[
F_2(\vphi_{\beta}, \vphi_A; n)\cap I_{n-1}(\bfa)
=
\bigcup_{1\leq a_n<\beta^n} \bigcup_{a_{n+1}\geq \frac{A^n}{a_{n} }} I_{n+1}(\bfa, a_{n}, a_{n+1}).
\]
Since
\[
\left| \bigcup_{a_{n+1}\geq \frac{\beta^n}{a_{n} }} I_{n+1}(\bfa, a_{n}, a_{n+1})\right|
\asymp 
\frac{1}{q_{n-1}^2(\bfa)a_nA^n},
\]
Lemma \ref{Lem:HWX:s:estimate} tells us that
\[
\sum_{1\leq a_n<\beta^n}
\frac{1 }{q_{n-1}^{2t}(\bfa)A^{nt} a_n^{t}}
\asymp 
\frac{(1-t)n\log \beta }{q_{n-1}^{2t}(\bfa)A^{nt} \beta^{n(t-1)}}
\quad\text{ for all } 0<t<1.
\]
Choosing $0<s<1$, $t=s$, and $\beta=A^s$, we obtain Lemma \ref{Lem:Bnd:HWX2} below.

\begin{lem01}\label{Lem:Bnd:HWX2}
For $A>1$, $0<s<1$ and $\veps>0$, there are $0<\veps'<\veps$ and $N\in\Na$ such that for all $n\in\Na_{\geq N}$ and $\bfa\in \Na^{n-1}$ there exists a countable cover $\mathfrak{C}(\bfa)$ of $E_2(\vphi_A;n)\cap I_{n-1}(\bfa)$ satisfying
\[
\Lambda_{s+\veps}(\mathfrak{C}(\bfa))\ll \frac{1}{q_{n-1}^{2(s+\veps')}(\bfa)A^{n(s^2+\veps')}}.
\]
\end{lem01}
We show the upper bound by working separately with the cases determined by the length of the overlaps of blocks with large products. 
In other words, we cover $\clF_3(\vphi_B)$ by
\begin{align*}
\clF_{3}^{0}
&:=
\left\{x\in [0,1): 
\begin{matrix}
a_{k}(x)a_{k+1}(x)a_{k+2}(x)\geq B^n, \\
a_{n}(x)a_{n+1}(x)a_{n+2}(x)\geq B^n,
\end{matrix}
\text{ for some } 1\leq k\leq n-3 \text{ for i.m.   } n\in\Na \right\}, \\
\clF_{3}^{1}
&:=
\left\{ x\in [0,1): 
\begin{matrix}
a_{n-2}(x) a_{n-1}(x) a_{n}(x) \geq B^n, \\
a_{n}(x) a_{n+1}(x) a_{n+2}(x) \geq B^n 
\end{matrix}
\text{ for i.m.   }n\in\Na\right\}, \\
\clF_{3}^{2}
&:=
\left\{ x\in [0,1): 
\begin{matrix}
a_{n-1}(x) a_{n}(x)a_{n+1}(x) \geq B^n, \\
a_{n}(x)a_{n+1}(x)a_{n+2}(x) \geq B^n 
\end{matrix}
\text{ for i.m.   } n\in\Na \right\}
\end{align*}
and use 
\begin{equation}\label{Eq:UB}
\hdim \clF_3(\vphi_B)
\leq 
\max\left\{
\hdim  \clF_3^2, \hdim  \clF_3^1, \hdim  \clF_3^0
\right\}.
\end{equation}

Define the functions $X,Y, Z:[0,1]\to[0,1]$ by 
\[
X(s) = \frac{s^2}{s^2-s+1},
\quad
Y(s) = \frac{s-s^2}{s^2-s+1},
\:\text{ and } \:
Z(s) = 1 - X(s) - Y(s),
\]
and the numbers 
\[
x_1=X(s_B), \quad
y_1=Y(s_B)
\quad\text{ and, } 
z_1=Z(s_B)
\]
(cfr. \eqref{Eq:DefXY}). 
Observe that, for $\frac{1}{2}<s<1$, 
\begin{equation}\label{Eq:EquivalenciaDeg}
g(s) =  2s-1 +Y(s)(1-s).
\end{equation}
By comparing the potentials in the formulas for $\hdim \clF_2(\vphi)$ and $\hdim \clE_{3}(\vphi)$ to $g$, we conclude that $\frac{1}{2}<s_B<1$. 
Several arguments are repeated multiple times. 
Then, aiming to keep the length of this paper reasonable, we only provide a full proof for $\clF_3^0$ and content ourselves with providing the essentials for the other sets.

\subsubsection{No overlaps}
We estimate $\hdim \clF_3^0(B)$.
For $n\in\Na$, we cover 
\[
\clF_{3,n}^{0}
=
\left\{ 
x\in [0,1): 
\begin{matrix}
\exists k\in\{1, \ldots, n-3\} \quad B^n\leq a_{k}(x)a_{k+1}(x)a_{k+2}(x) ,\\
B^n\leq a_{n}(x)a_{n+1}(x)a_{n+2}(x)
\end{matrix}
\right\}
\]
by the sets 
\[
A_{n}^{0}(r)
=
\left\{ 
x\in [0,1):
\begin{matrix}
\exists k\in\{1, \ldots, n-3\}\quad B^n\leq a_{k}(x)a_{k+1}(x)a_{k+2}(x),  \\
a_n(x)a_{n+1}(x)\leq B^{nr}, \; B^n\leq a_n(x)a_{n+1}(x)a_{n+2}(x) 
\end{matrix}
\right\},
\]
\[
A_{n}^{1}(r)
=
\left\{ 
x\in [0,1):
\begin{matrix}
\exists k\in\{1, \ldots, n-3\}\quad B^n\leq a_{k}(x)a_{k+1}(x)a_{k+2}(x),  \\
B^{nr}\leq a_n(x)a_{n+1}(x) 
\end{matrix}
\right\},
\]
where $0<r<1$ is to be determined; therefore, 
\begin{equation}\label{Eq:UB:F03}
\clF_3^0 
=
\limsup_{n\to\infty}\clF_{3,n}^{0}
\subseteq 
\limsup_{n\to\infty} A_{n}^{0}(r)
\cup
\limsup_{n\to\infty}  A_{n}^{1}(r).
\end{equation}
let $s$ be such that $\frac{1}{2}<s_{B}<s<1$, and write $r=\frac{s}{s^2-s+1}$.
Direct computations yield
\begin{equation}\label{Eq:ol0:01}
F(s)=
2s-1+rs-r+s
=
2 s - 1 + \frac{s^2 - s}{s^2-s+1} + \frac{s^3-s^2+s}{s^2-s+1}
>
g(s).
\end{equation}
For each $1\leq k\leq n-3$, define 
\[
A_{n,k}^{0}
=
\left\{ 
x\in [0,1): 
\begin{matrix}
B^n\leq a_{n}(x)a_{n+1}(x)a_{n+2}(x) ,\, a_n(x)a_{n+1}(x)\leq  B^{nr},\\
B^n\leq a_{k}(x)a_{k+1}(x)a_{k+2}(x) 
\end{matrix}
\right\}, 
\]
so
\[
A_{n}^{0} = \bigcup_{k=1}^{n-3} A_{n,k}^{0}.
\]
For $1\leq k\leq n-3$ and $\bfa=(a_1,\ldots, a_{n+1})\in \Na^{n+1}$ such that $a_ka_{k+1}a_{k+2}\geq B^n$ and $a_na_{n+1}\leq B^{nr}$, define 
\[
J_{n+1}(\bfa) = \bigcup_{a_{n+2}\geq \frac{B^{n}}{a_na_{n+1}}} I_{n+2}(\bfa a_{n+2}),
\quad\text{ so }\quad
\left| J_{n+1}(\bfa)\right| \asymp \frac{1}{q_{n-1}^2(a_1,\ldots, a_{n-1})a_na_{n+1}B^{n}}.
\]
Then, $\mathfrak{A}_{n,k}^{0} = \left\{ J_{n+1}(\bfa): a_ka_{k+1}a_{k+2}\geq B^n, \,a_na_{n+1}\leq B^{nr}\right\}$ is a covering of $A_{n,k}^{0}$.
Next, we estimate $\Lambda_s(\mathfrak{A}_{n,k}^{0})$.
First, for  $\bfa=(a_{1}, \ldots, a_{n-1})\in\Na^{n-1}$ such that $a_ka_{k+1}a_{k+2}\geq B^n$, we apply Lemma \ref{Lem:HWX:s:estimate} to get
\begin{align*}
\sum_{a_na_{n+1}\leq B^{nr}} \left| J_{n+1}(\bfa)\right|^{s}
&\asymp
\sum_{a_na_{n+1}\leq B^{nr}} \frac{1}{q_{n-1}^{2s}(a_1, \ldots, a_{n-1}) a_n^sa_{n+1}^sB^{ns}} \\
&\asymp
\frac{nr\log B}{q_{n-1}^{2s}(a_1, \ldots, a_{n-1})  B^{nr(s-1)}B^{ns}}.
\end{align*}
Second, we allow $(a_k, a_{k+1}, a_{k+2})$ to vary along the triples such that $a_ka_{k+1}a_{k+2}\geq B^n$ leaving all the other variables constant:
\[
\sum_{a_k a_{k+1} a_{k+2}\geq B^n}
\sum_{a_na_{n+1}\leq B^{nr}} 
\left| J_{n+1}(\bfa)\right|^{s}
\asymp 
\frac{n^2r\log B}{q_{n-4}^{2s}(a_1, \ldots, a_{k-1}, a_{k+4}, \ldots, a_{n-1}) B^{2s-1}  B^{n(rs-r+s)}}.
\]
Lastly, if $\bfb=(a_1, \ldots, a_{k-1}, a_{k+3}, \ldots, a_{n-1})$ runs along $\Na^{n-4}$, we obtain
\begin{align}
\Lambda_s(\mathfrak{A}_{n,k}^0)
&\asymp 
\sum_{\bfb\in\Na^{n-4}}
\sum_{a_k a_{k+1} a_{k+2}\geq B^n}
\sum_{a_na_{n+1}\leq B^{nr}} 
\left| J_{n+1}(\bfa)\right|^{s} \nonumber \\
&\asymp 
\sum_{\bfb\in\Na^{n-4}}
\frac{n^2r\log B}{q_{n-4}^{2s}(\bfb ) B^{n(2s-1+rs-r+s)}}. \label{Eq:UP:Ref:01}
\end{align}
Put $\veps=s-s_B>0$ and define $\tau>0$ by $F(s_B)=g(s_B)+\tau$ with $F$ as in \eqref{Eq:ol0:01}. 
Because of \eqref{Eq:DefSmB}, \eqref{Eq:ol0:01}, and the continuity and strict growth of $g$ and $F$, large values of $n$ satisfy $s_{n,B}<s_B+\frac{\veps}{2}$ and 
\[
F(s)
>
F\left(s_{n,B} + \frac{\veps}{2}\right)
> 
g(s_{n,B})+\frac{\tau}{2}.
\]
Large values of $n$ also satisfy $nr\log B\leq  B^{\frac{n\tau}{2}}$.
Then, writing $\delta=\frac{1}{4}\min\{\tau, \veps\}>0$, 
\begin{align*}
\Lambda_s(\mathfrak{A}_{n,k}^0)
&\ll 
\sum_{\bfb\in\Na^{n-4}}
\frac{1}{q_{n-4}^{2(s_{n,B}+\frac{\veps}{2})}(\bfb ) B^{n(g(s_{n,B}) + \frac{\tau}{2} ) } } \\
&\leq 
B^4
\sum_{\bfb\in\Na^{n-4}}
\frac{1}{q_{n-4}^{2(s_{n,B}+\frac{\veps}{2})}(\bfb ) B^{(n-4)(g(s_{n,B})+\frac{\tau}{2})} }  \\
&\leq 
B^4\frac{1}{B^{(n-4) \delta  }} \\
&\ll_{B} 
\frac{1}{B^{\delta n }}
\end{align*}
For all $N\in\Na$, the collection
\[
\mathfrak{A}_{N}^{0} = \bigcup_{n\geq N} \bigcup_{k=1}^{n-3} \mathfrak{A}_{n,k}^0.
\]
is a countable covering of $\displaystyle\limsup_{n\to\infty} A_n^{0}$, and, if $N$ is large enough, 
\begin{align*}
\Lambda_s(\mathfrak{A}_{N}^{0})
&=
\sum_{n= N}^{\infty}
\sum_{k=1}^{n-3} 
\Lambda_s\left( \mathfrak{A}_{n,k}^0\right) \\
&\ll 
\sum_{n= N}^{\infty}
\frac{1}{ B^{\delta n }} = \frac{1}{B^{N\delta}} \left(\frac{B^{\delta}}{1-B^{\delta}}\right).
\end{align*}
 Since $s>s_B$ was arbitrary, $\hdim \displaystyle\limsup_{n\to\infty} A_n^{0}\leq s_B$.

Next, we obtain an upper bound for the Hausdorff dimension of the limsup of $(A_n^1)_{n\geq 1}$.
Take real numbers $s$, $\veps$, and $\tilde{s}$ such that
\[
s_B<s<1, \quad 
\veps=s-s_B>0, \quad 
\text{ and }\quad 
s_B<\tilde{s}<s.
\]
Given $n\in\Na_{\geq 4}$, cover $A_{n}^{1}$ by the sets
\[
A_{n,k}^{1}=
\left\{x\in [0,1): a_{k}(x)a_{k+1}(x)a_{k+2}(x)\geq B^n, \; a_{n}(x)a_{n+1}(x)\geq B^{nr}\right\}, 
\;
1\leq k\leq n-3.
\]
For such $k$ and $\bfa=(a_1,\ldots, a_{n-1})\in\Na^{n-1}$ satisfying $a_ka_{k+1}a_{k+2}\geq B^n$, let $\mathfrak{A}_{n,k}^1(\bfa) $ be the covering of 
\[
\bigcup_{\substack{a_n, a_{n+1} \in\Na \\ a_{n}a_{n+1}\geq B^{rn} }}
I_{n+1}(\bfa, a_n, a_{n+1})
\]
obtained for $\tilde{s}$ from Lemma \ref{Lem:Bnd:HWX2}.
This means that there is some $0<\veps'<\veps$ such that, if $n$ is large,
\[
\Lambda_s\left( \mathfrak{A}_{n,k}^1(\bfa)\right)
\ll 
\frac{1}{q_{n-1}^{2s}(\bfa)B^{nr(\tilde{s}^2+\veps')}}.
\]
Writing $\bfb=(a_1, \ldots, a_{k-1},a_{k+3}, \ldots, a_{n-1})$ and using Lemma \ref{Lem:s:estimate:02}, we get
\begin{align}
\sum_{a_ka_{k+1}a_{k+2}\geq B^n}
\Lambda_s\left( \mathfrak{A}_{n,k}^1(a_1, \ldots, a_{n-1})\right)
&\ll 
\frac{n^2\log^2 B}{q_{n-4}^{2s}(\bfb)B^{n(2s-1)}B^{nr(\tilde{s}^2+\veps')}} \nonumber\\
&<
\frac{n^2\log^2 B}{q_{n-4}^{2s}(\bfb)B^{n(2\tilde{s}-1)+nr\tilde{s}^2} } \nonumber\\
&<
\frac{n^2\log^2 B}{q_{n-4}^{2s}(\bfb)B^{ng(\tilde{s})}}. \label{Eq:UP:Ref:02}
\end{align}
Since $g$ is strictly increasing and $s_B<\tilde{s}$, we may define $\delta'>0$ by $g(\tilde{s})=g(s_B)+3\delta'$; hence, $g(\tilde{s}) > g(s_{n,B}) + 2\delta'$ for large $n\in\Na$. 
Large values of $n$ also satisfy $n^2\log^2B< B^{n\delta'}$, so
\[
\sum_{\bfb\in\Na^{n-4}}
\sum_{a_ka_{k+1}a_{k+2}\geq B^n}
\Lambda_s\left( \mathfrak{A}_{n,k}^1(a_1, \ldots, a_{n-1})\right)
\ll
\frac{1}{B^{\delta'n}}.
\]
We define
\begin{align*}
\mathfrak{A}_{n,k}^{1}
&=
\left\{ \mathfrak{A}_{n,k}^1(\bfa): \bfa\in\Na^{n-1}, \, a_{k}a_{k+1}a_{k+2}\geq B^n\right\}, \;\text{ and }\;\\    
\mathfrak{A}_{N}^{1} 
&=
\bigcup_{n\geq N}
\bigcup_{k=1}^{n-3}
\mathfrak{A}_{n,k}^{1} \;\text{ for } N\in\Na_{\geq 4}.
\end{align*}
Finally, for large $N$, the collection $\mathfrak{A}_{N}^{1}$ is a countable cover of $\displaystyle\limsup_{n\to\infty} A_n^{1}$ such that $\Lambda_s(\mathfrak{A}_{N}^{1})\ll_{\delta'} B^{-N\delta'}$, and, thus, $\hdim\displaystyle\limsup_{n\to\infty} A_n^{1}\leq s_B$.
Therefore, by virtue of \eqref{Eq:UB:F03}, $\hdim \clF_3^0\leq s_B$.
\subsubsection{Overlaps of length $1$}
The keys to the previous case were inequalities \eqref{Eq:UP:Ref:01} and \eqref{Eq:UP:Ref:02}.
Once we got them, we used the convergence of $(s_{n,B})_{n\geq 1}$, and the continuity and strict growth of the functions involved. 
These last steps are straightforward; the main complication is constructing the adequate coverings of $\clF_0^3$, $\displaystyle\limsup_{n\to\infty} A_n^0$, and $\displaystyle\limsup_{n\to\infty} A_n^1$. 
With this in mind, we can significantly trim the forthcoming proofs without jeopardizing clarity nor validity.
For $n\in\Na$, define
\begin{align*}
V_{n}^{0} &=
\left\{ x\in [0,1): B^{n(1-y_1)}\leq a_{n-2}(x)a_{n-1}(x), \,  B^{ny_1}\leq a_n(x) \right\}, \\
V_{n}^{1}&=
\left\{ x\in [0,1): B^{n(1-y_1)}\leq a_{n-2}(x)a_{n-1}(x), \, a_n(x)< B^{ny_1}, \, B^n\leq a_{n}(x)a_{n+1}(x)a_{n+2}(x) \right\},   \\
V_{n}^{2}&=
\left\{ x\in [0,1): a_{n-2}(x)a_{n-1}(x)\leq B^{n(1-y_1)}, \, B^n\leq a_{n-2}(x)a_{n-1}(x) a_n(x) \right\},
\end{align*}
then
\begin{equation}\label{Eq:UB:F13}
\clF_{3}^{1} 
= 
\limsup_{n\to\infty} V_{n}^{0}
\cup
\limsup_{n\to\infty} V_{n}^{1}
\cup
\limsup_{n\to\infty} V_{n}^{2}.
\end{equation}

To bound the Hausdorff dimension of the limsup of $(V_{n}^{0})_{n\geq 1}$, for each $\bfa\in\Na^{n-3}$ define
\[
J_{n-1}(\bfa a_{n-2}a_{n-1})
=
\bigcup_{a_n\geq B^{ny_1}} I_n(\bfa a_{n-2}a_{n-1}a_n),
\;\text{ so }\;
| J_{n-1}(\bfa a_{n-2}a_{n-1})|
\asymp 
\frac{1}{q_{n-3}^2(\bfa)a_{n-2}^{2}a_{n-1}^{2}B^{ny_1}},
\]
and, for $\frac{1}{2}<s<1$, 
\begin{align*}
\sum_{a_{n-2} a_{n-1} \geq B^{n(1-y_1)}} | J_{n-1}(\bfa a_{n-2}a_{n-1})|^s
&\asymp 
\sum_{a_{n-2} a_{n-1}\geq B^{n(1-y_1)}} \frac{1}{q_{n-3}^{2s}(\bfa)a_{n-2}^{2s}a_{n-1}^{2s}B^{nsy_1}} \\
&\asymp
\frac{1}{q_{n-3}^{2s}(\bfa)B^{n(1-y_1)(2s-1)}B^{nsy_1}} \\
&\asymp
\frac{1}{q_{n-3}^{2s}(\bfa)B^{n( (2s-1) +y_1(1-s)) }}.
\end{align*}
Using the strict growth of $g$ and \eqref{Eq:EquivalenciaDeg}, we conclude $\displaystyle\limsup_{n\to\infty} V_{n}^{0}\leq s_B$.\par

We now deal with $(V_{n}^{1})_{n\geq 1}$.
For $\bfa=(a_1, \ldots, a_{n-3}) \in \Na^{n-3}$, $a_{n-2}$ and $a_{n-1}$ such that $a_{n-2}a_{n-1}\geq B^{n(1-y_1)}$ and $a_n<B^{ny_1}$, Lemma \ref{Lem:Bnd:HWX2} gives us a covering $\mathfrak{D}(\bfa, a_{n-2}, a_{n-1}, a_n)$ of 
\[
\bigcup_{a_{n+1}a_{n+2}\geq \frac{B^n}{a_n}} I_{n+2}(a_1, \ldots, a_{n+2}) 
\]
such that, for a fixed $\tilde{s}$ with $s_B<\tilde{s}<s$, 
\[
\Lambda_s(\mathfrak{D}(\bfa, a_{n-2}, a_{n-1}, a_n))
\ll
\frac{1}{q_{n-3}^{2\tilde{s}}(a_1, \ldots, a_{n-3}) a_{n-2}^{2\tilde{s}}a_{n-1}^{2\tilde{s}}a_{n}^{2\tilde{s}-\tilde{s}^2}B^{n\tilde{s}^2}}.
\]
Since $0<2\tilde{s}-\tilde{s}^2<1$, we have that
\begin{align*}
\sum_{a_n \leq  B^{ny_1}} \Lambda_s(\mathfrak{D}(\bfa, a_{n-2}, a_{n-1}, a_n))
&\asymp 
\frac{1}{q_{n-3}^{2\tilde{s}}(\bfa) a_{n-2}^{2\tilde{s}}a_{n-1}^{2\tilde{s}}B^{ny_1(2\tilde{s}-\tilde{s}^2-1)}B^{n\tilde{s}^2}} \\
&=
\frac{1}{q_{n-3}^{2\tilde{s}}(\bfa) a_{n-2}^{2\tilde{s}}a_{n-1}^{2\tilde{s}}B^{n(-y_1(1-\tilde{s})^2+\tilde{s}^2)}}.
\end{align*}
Then, by Lemma \ref{Lem:HWX:s:estimate}, 
\begin{align*}
\sum_{a_{n-2}a_{n-1}\geq B^{n(1-y_1)}}
\sum_{a_n \leq  B^{ny_1}} \Lambda_s(\mathfrak{D}(\bfa, a_{n-2}, a_{n-1}, a_n))
&\ll 
\sum_{a_{n-2}a_{n-1}\geq B^{n(1-y_1)}}
\frac{1}{q_{n-3}^{2\tilde{s}}(\bfa) a_{n-2}^{2\tilde{s}}a_{n-1}^{2\tilde{s}}B^{n(-y_1(1-\tilde{s})^2+\tilde{s}^2)}} \\
&\asymp 
\frac{1}{q_{n-3}^{2\tilde{s}}(\bfa) B^{n(2\tilde{s}-1)(1-y_1)}B^{n(-y_1(1-\tilde{s})^2+\tilde{s}^2)}} \\
&=
\frac{1}{q_{n-3}^{2\tilde{s}}(\bfa) B^{n(2\tilde{s}-1 +\tilde{s}^2(1-y_1))}} \\
&<
\frac{1}{q_{n-3}^{2\tilde{s}}(\bfa) B^{ng(\tilde{s})}}
\end{align*}
and we can conclude that $\hdim\displaystyle\limsup_{n\to\infty} V_{n}^{1}\leq s_B$.

Next, we focus on $(V_{n}^{2} )_{n\geq 1}$. 
For each $\bfa\in \Na^{n-1}$ such that $a_{n-2}a_{n-1}\leq B^{n(1-y_1)}$ define 
\[
J_n(\bfa)
=
\bigcup_{a_n\geq \frac{B^{n }}{a_{n-2}a_{n-1} }} I_{n}(\bfa, a_n),
\;\text{ so }\;
\left| J_n(\bfa)\right|
\asymp
\frac{1}{q_{n-3}^2(a_1, \ldots, a_{n-3}) a_{n-2}a_{n-1}B^{n}}.
\]
From Lemma \ref{Lem:HWX:s:estimate}, we obtain 
\begin{align*}
\sum_{a_{n-1}a_{n-2}\leq B^{n(1-y_1)}}
|J_n(a_1,\ldots, a_{n-1})|^s
&\asymp 
\frac{n(1-y_1)\log B}{q_{n-3}^{2s}(a_1, \ldots, a_{n-3}) B^{ns}B^{n(1-y_1)(1-s)}} \\
&=
\frac{n(1-y_1)\log B}{q_{n-3}^{2s}(a_1, \ldots, a_{n-3}) B^{n(2s-1 + y_1(1-s)) }} \\
&=
\frac{n(1-y_1)\log B}{q_{n-3}^{2s}(a_1, \ldots, a_{n-3}) B^{ng(s)}}.
\end{align*}
Therefore, $\hdim \displaystyle\limsup_{n\to\infty} V_{n}^{2}\leq s_B$ and, by \eqref{Eq:UB:F13}, $\hdim \clF_3^1\leq s_B$.
 \subsubsection{Overlaps of length 2}
To construct a useful covering of $\clF_{3}^{2}$, for every $n\in\Na$ define
\begin{align*}
C_{0,n}&=
\left\{ x\in [0,1):  B^{n}\leq a_{n-1}(x) \right\}, \\
C_{1,n}&=
\left\{ x\in [0,1): a_{n-1}(x)\leq B^{nz_1}, \, B^n\leq a_{n}(x)a_{n+1}(x) a_{n+2}(x) \right\}, \\
C_{2,n}&=
\left\{ x\in [0,1): B^{nz_1}\leq a_{n-1}(x)\leq B^{n}, \,  B^n \leq a_{n-1}(x)a_{n}(x) a_{n+1}(x), B^{nz_1}\leq a_{n+2}(x) \right\}, \\
C_{3,n}&=
\left\{ x\in [0,1): B^{nz_1}\leq a_{n-1}(x)\leq B^{n}, \, B^n \leq a_n(x)a_{n+1}(x) a_{n+2}(x), a_{n+2}(x)\leq B^{nz_1} \right\}; 
\end{align*}
hence, 
\begin{equation}\label{Eq:UB:F23}
\clF_{3}^{2}
= 
\limsup_{n\to\infty} C_{0,n}
\cup 
\limsup_{n\to\infty} C_{1,n}
\cup 
\limsup_{n\to\infty} C_{2,n}
\cup 
\limsup_{n\to\infty} C_{3,n}.
\end{equation}
Since $g(s)<s$ when $\frac{1}{2}<s<1$, Theorem \ref{TEO:HWX:HD} implies $\displaystyle\hdim \limsup_{n\to\infty } C_{0,n} \leq  s_B$.\par 

Let us discuss the limsup of $(C_{1,n})_{n\geq 1}$.
Let $f_3:[0,1]\to\RE$, $f_3(s) = \frac{s^3}{s^2-s+1}$, be as in Theorem \ref{TEO:HWX:HD}.
Direct computations show that
\[
g(s) < (2s-1)Z(s) + f_3(s) \; \text{ for }\; \frac{1}{2}<s<1.
\]
Let $s, \tilde{s}$ be such that $s_B<\tilde{s}<s<1$.
As we did in Lemma \ref{Lem:Bnd:HWX2}, we can follow the proof of \cite[Theorem 1.7]{HuangWuXu2020} to construct for all large $n\in\Na$ and any $\bfa\in \Na^{n-2}$, $a_{n-1} \in\Na$ a cover $\mathfrak{G}(\bfa, a_{n-1}, a_n)$ of
\[
\bigcup_{a_na_{n+1}a_{n+2}\geq B^{n}} I_{n+2}(\bfa, a_{n-1}, a_{n}, a_{n+1}, a_{n+2})
\]
such that 
\[
\Lambda_s(\mathfrak{G}(\bfa, a_{n-1}))
\ll 
\frac{1}{q_{n-2}^{2\tilde{s}}(\bfa)a_{n-1}^{2\tilde{s} } B^{nf_3(\tilde{s}) }};
\]
hence, assuming $s$ is close to $s_B$, 
\begin{align*}
\sum_{a_{n-1}\leq B^{nz_1}}
\Lambda_s(\mathfrak{G}(\bfa, a_{n-1}))
&\ll 
\frac{1}{q_{n-2}^{2\tilde{s}}(\bfa)B^{nz_1(2\tilde{s}-1) } B^{nf_3(\tilde{s}) }} \\
&=
\frac{1}{q_{n-2}^{2\tilde{s}}(\bfa) B^{n(z_1(2\tilde{s}-1) + f_3(\tilde{s})) }} 
<
\frac{1}{q_{n-2}^{2\tilde{s}}(\bfa) B^{ng(\tilde{s})}}.
\end{align*}
The estimate $\hdim \displaystyle\limsup_{n\to\infty} C_{2,n}\leq s_B$ follows from $g(s_B) = 2 s_B-1+s_Bz_1$, the strict growth of $s\mapsto 2s-1+sz_1$, and
\begin{align*}
\sum_{B^{nz_1}\leq a_{n-1}\leq B^n}
\sum_{\frac{B^{n}}{a_{n-1}}\leq a_{n}a_{n+1}}
\frac{1}{a_{n-1}^{2s}a_{n}^{2s}a_{n+1}^{2s}B^{nsz_1}}
&\ll
\sum_{B^{nz_1}\leq a_{n-1}\leq B^n}
\frac{1}{a_{n-1}^{2s}\left(\frac{B^n}{a_{n-1}}\right)^{2s-1} B^{nsz_1}} \\
&=
\sum_{B^{nz_1}\leq a_{n-1}\leq B^n}
\frac{1}{a_{n-1}  B^{n(2s-1 +sz_1)}} \\
&\leq
\frac{n\log B}{B^{n(2s-1 +sz_1)}}.
\end{align*}
Lastly, note that $n\in\Na$,
\[
C_{3,n}
\subseteq 
C_{3,n}'
=
\left\{ x\in [0,1): B^{nz}\leq a_{n-1}(x)\leq B^{n}, \, B^{n(1-z_1)} \leq a_{n}(x) a_{n+1}(x) \right\}.
\]
By Lemma \ref{Lem:HWX:s:estimate}, for each $\bfa\in\Na^{n-2}$ and $B^{nz_1}\leq a_{n-1}(x)\leq B^{n}$ there is a countable cover $\mathfrak{C}'(\bfa, a_{n-1})$  of $C_{3,n}'\cap I_{n-1}(\bfa, a_{n-1})$ such that
\begin{align*}
\sum_{B^{nz_1}\leq a_{n-1}(x)\leq B^{n}}
\sum_{C\in \mathfrak{C}'(\bfa, a_{n-1})}
|C|^s
&\ll
\sum_{B^{n z_1}\leq a_{n-1}(x)\leq B^{n}}
\frac{1}{q_{n-2}^{2s}(\bfa)a_{n-1}^{2s}B^{ns^2(1-z_1)}}\\
&\leq
\sum_{B^{nz}\leq a_{n-1}(x)}
\frac{1}{q_{n-2}^{2s}(\bfa)a_{n-1}^{2s}B^{ns^2(1-z_1)}}\\
&\leq
\frac{1}{q_{n-2}^{2s}(\bfa)B^{n(2s-1)z}B^{ns^2(1-z_1)}}\\
&=
\frac{1}{q_{n-2}^{2s}(\bfa)B^{n(z_1(2s-1) + s^2(1-z_1)) }}.
\end{align*}
Because $s\mapsto z_1(2s-1) + s^2(1-z_1)$ is strictly increasing and $g(s_B)=z_1(2s_B-1) + s_B^2(1-z_1)$, we conclude
\[
\hdim\limsup_{n\to\infty} C_{3,n}
\leq 
\hdim\limsup_{n\to\infty} C_{3,n}'
\leq 
s_B.
\]
Therefore, \eqref{Eq:UB:F23} yields $\hdim \clF_3^2\leq s_B$.

\section{Final remarks and further results}\label{Section:FinalRemarks}

In the introduction, we said that Theorem \ref{TEO:MAIN:MEASURE} is an important piece in the extension of the laws of large numbers \eqref{Eq:Thm:DV} and \eqref{Eq:Thm:HHY}. 
In a forthcoming article, we study the asymptotic behavior of the sums of products of consecutive partial quotients. 
More precisely, for $n,\ell\in\Na$ and $x\in [0,1)\setminus\QU$, define
\[
S_{n,\ell}(x)=\sum_{j=1}^n a_{j}(x) \cdots a_{j+\ell-1}(x).
\]
\begin{conj01}\label{TEO:PROMISED:01}\
For every $\veps>0$, we have
\[
\lim_{n\to \infty}
\leb
\left(\left\{ 
x\in [0,1): 
\left| \frac{1}{n\log^{\ell}n} S_{n,\ell}(x) - \frac{1}{\ell \log 2}\right|\geq \veps\right\}\right)=0.
\]
Also, almost every $x\in [0,1)$ satisfies
\[
\lim_{n\to\infty} 
\frac{1}{n\log^{\ell} n} 
\left( S_{n,\ell}(x) - \max_{1\leq j \leq n} a_{j}(x)\cdots a_{j+\ell-1}(x)\right)
=
\frac{1}{\ell \log 2}.
\]
\end{conj01}

Regarding the Hausdorff dimension of certain Lebesgue-null sets, for a non-decreasing function $\vphi:\Na\to\RE$ such that $\vphi(n)\to\infty$ as $n\to\infty$, we consider
\[
G_{\ell}(\vphi)
:=
\left\{ x\in [0,1): \lim_{n\to\infty} \frac{S_{n,\ell}(x)}{\vphi(n)}=1\right\}.
\]
\begin{conj01}
Let $\alpha$ be a positive real number.
\begin{enumerate}[\rm (I)]
    \item If $0<\alpha<\frac{1}{2}$ and $\vphi(n)=e^{n^{\alpha}}$, then $\hdim G_{\ell}(\vphi)=1$.
    \item If $\frac{1}{2}\leq \alpha$ and $\vphi(n)=e^{n^{\alpha}}$, then $\hdim G_{\ell}(\vphi)=\frac{1}{2}$.
    \item If $1<\alpha$ and $\vphi(n)=e^{\alpha^n}$, then $\hdim G_{\ell}(\vphi)=\frac{1}{1+\alpha}$.
\end{enumerate}
\end{conj01}
 Recently, Song, Yu, and Yu \cite{SYY2025} investigated the Hausdorff dimension of a related set concerning the product of two partial quotients. 
 To define it, for $x\in [0,1)\setminus\QU$ and $n\in\Na$, put
 \[
 L_n(x) = \max_{1\leq i\leq n} \, a_i(x)a_{i+1}(x),
 \]
  For a non-decreasing function $\vphi:\RE_{>0}\to\RE_{>0}$ satisfying certain regularity properties, they computed the Hausdorff dimension of 
        \[
            L(\varphi):= \left\{ x\in(0,1): \lim_{n\to \infty} \frac{L_n(x)}{\varphi(n)}=1\right\}.
        \] 
        Their work leads us to the following problem.
        \begin{prob01}
        Let $\ell\in\Na_{\geq 3}$.
        Can we extend the results in \cite{SYY2025} to the product of $\ell$ consecutive partial quotients? 
        That is, for $x\in [0,1)\setminus\QU$ and $n\in\Na$, define
        \[
        L_{\ell,n}(x) = \max_{1\leq i\leq n}\, a_i(x)a_{i+1}(x)\cdots a_{i+ \ell -1 }(x) .
        \]
        Given $\rho>0$ and a function $\vphi:\RE_{>0}\to\RE_{>0}$ such that $\log \vphi$ is regularly increasing of index $\rho$, determine the Hausdorff dimension of the set
        \[
            L_{\ell}(\varphi)
            := 
            \left\{ x\in(0,1): \lim_{n\to \infty} \frac{L_{\ell,n}(x)}{\varphi(n)}=1\right\}.
        \] 
        \end{prob01}
We conclude our work with two remarks concerning our main theorems.
\begin{rema01}
The $n$th term of one of the series in Theorem \ref{TEO:MAIN:MEASURE} is
\[
\frac{n\log^{2(\ell-1)}\vphi(n)}{\vphi^2(n)}
+
\frac{\log^{\ell-2} \vphi(n)}{\vphi(n)}.
\]
The first term accounts for large products that appear in non-overlapping blocks, and the second term corresponds to overlapping blocks. Hence, for $\ell=1$ the second term is non-existent. 
\end{rema01}

\begin{rema01}
It seems that the proof of Theorem \ref{TEO:MAIN:HD} can be adapted to compute $\hdim \clF_{\ell}(\vphi)$ for $\ell\geq 4$.
A recursive formula for the potential giving $\hdim \clF_{\ell}(\vphi)$ could be expected, we have not used the formulas corresponding to $\hdim \clF_1(\vphi_B)$ or $\hdim \clF_2(\vphi_B)$. Instead, we used those for $\hdim \clE_1(\vphi_B)$ and $\hdim \clE_2(\vphi_B)$.
\end{rema01}

\bibliography{References}
\bibliographystyle{abbrv}
\end{document}